\documentclass{gtart_h}

\def\ifplaintex{\expandafter\ifx\csname documentclass\endcsname\relax}

\def\gtp{{\mathsurround=0pt\it $\cal G\mskip-2mu$eometry \&\ 
$\cal T\!\!$opology $\cal P\!$ublications}}  

\def\recd{{\small Received:\qua\receiveddate\ifx\reviseddate\relax
\else\qquad Revised:\qua\reviseddate\fi\par}} 

\def\lognumber#1{\def\thelognumber{#1}}
\def\volumenumber#1{\def\thevolumenumber{#1}}
\def\volumeyear#1{\def\thevolumeyear{#1}}
\def\papernumber#1{\def\thepapernumber{#1}}
\def\pagenumbers#1#2{\def\startpage{#1}\def\finishpage{#2}}
\def\published#1{\def\publishdate{#1}}

\def\received#1{\def\receiveddate{#1}}
\def\revised#1{\def\reviseddate{#1}}
\def\accepted#1{\def\accepteddate{#1}}
\def\asciititle#1{\def\theasciititle{#1}}

\let\\\par\let\thelognumber\relax\let\thevolumenumber\relax
\let\thepapernumber\relax\let\thevolumeyear\relax\let\startpage\relax
\let\finishpage\relax\let\publishdate\relax\let\receiveddate\relax
\let\reviseddate\relax\let\accepteddate\relax\let\theasciititle\relax
\let\theasciiauthors\relax
\let\theasciiabstract\relax

\let\theasciiemail\relax

\ifplaintex
\font\logobig=cmssbx10 scaled 3836
\font\logomed=cmssbx10 scaled 2557
\else
\font\logobig=cmssbx10 scaled 4200
\font\logomed=cmssbx10 scaled 2800
\fi

\long\def\makeagttitle{   
\count0=\startpage
\agt\hfill      
\hbox to 45truept{\vbox to 0pt{\vglue -13truept{\logomed A\kern -.37em{\logobig 
T}\kern -.38em G}\vss}\hss}
\break
{\small Volume \thevolumenumber\ (\thevolumeyear)
\startpage--\finishpage\nl
Published: \publishdate}

\vglue .25truein

{\parskip=0pt\leftskip 0pt plus
1fil\def\\{\par\smallskip}{\Large\bf\thetitle}\par\medskip} \vglue
0.05truein

{\parskip=0pt\leftskip 0pt plus 1fil\def\\{\par}{\sc\theauthors}
\par\medskip}%
 
\vglue 0.03truein 

{\small\leftskip 25truept\rightskip 25truept{\bf Abstract}\stdspace\theabstract

{\bf AMS Classification}\stdspace\theprimaryclass
\ifx\thesecondaryclass\relax\else; \thesecondaryclass\fi\par
{\bf Keywords}\stdspace \thekeywords\par}\vglue 7truept

}   

\ifplaintex
\font\phead=cmsl9 scaled 950
\font\pnum=cmbx10 scaled 913
\font\pfoot=cmsl9 scaled 950
\headline{\vbox to 0pt{\vskip -4.5mm\line{\small\phead\ifnum
\count0=\startpage ISSN 1472-2739 (on-line) 1472-2747 (printed)
\hfill {\pnum\folio}\else\ifodd\count0\def\\{ }%
\ifx\theshorttitle\relax\thetitle\else\theshorttitle\fi\hfill{\pnum\folio}
\else\def\\{ and }{\pnum\folio}\hfill\ifx\theshortauthors\relax\theauthors
\else\theshortauthors\fi\fi\fi}\vss}}
\footline{\vbox to 0pt{\vglue 0mm\line{\small\pfoot\ifnum\count0=\startpage
\copyright\ \gtp\hfill\else
\agt, Volume \thevolumenumber\ (\thevolumeyear)\hfill\fi}\vss}}
\else
\font=cmsl9 scaled 1050
\font\lnum=cmbx10 
\font=cmsl9 scaled 1050
\makeatletter
\def\@oddhead{{\small\lhead\ifnum\count0=\startpage ISSN 1472-2739 
(on-line) 1472-2747 (printed)\hfill {\lnum\number\count0}\else\ifodd\count0
\def\\{ }\ifx\theshorttitle\relax \thetitle \else\theshorttitle\fi\hfill
{\lnum\number\count0}\else\def\\{ and }{\lnum\number\count0}
\hfill\ifx\theshortauthors\relax 
\theauthors\else\theshortauthors\fi\fi\fi}}\def\@evenhead{\@oddhead}
\def\@oddfoot{\small\lfoot\ifnum\count0=\startpage\copyright\ \gtp\hfill\else
\agt, Volume \thevolumenumber\ (\thevolumeyear)\hfill\fi}
\def\@evenfoot{\@oddfoot}
\makeatother
\fi
\let\maketitlepage\makeagttitle
\let\makeshorttitle\maketitlepage
\let\maketitle\maketitlepage

\newwrite\gtoutfile
\long\gdef\makeheadfile{  
{\def\\{, }\def\s{ }
\immediate\openout\gtoutfile head.xxx
\immediate\write\gtoutfile{Proxy-for: \ifx\theasciiauthors\relax
\theauthors\else\theasciiauthors\fi\s<\ifx\theasciiemail\relax\theemail\else\theasciiemail\fi>}
\immediate\write\gtoutfile{\noexpand\\}
\immediate\write\gtoutfile{Authors: \ifx\theasciiauthors\relax
\theauthors\else\theasciiauthors\fi}
{\def\\{ }\immediate\write\gtoutfile{Title: \ifx\theasciititle\relax
\thetitle\else\theasciititle\fi}}
\immediate\write\gtoutfile{Subj-class: GT or SG, GR etc}
\immediate\write\gtoutfile{MSC-class: \theprimaryclass\ifx\thesecondaryclass\relax\else, \thesecondaryclass\fi}
\immediate\write\gtoutfile{Journal-ref: Algebr. Geom. Topol. \thevolumenumber\s
(\thevolumeyear) \startpage-\finishpage}
\immediate\write\gtoutfile{Comments: Published by Algebraic and
Geometric Topology at}
\immediate\write\gtoutfile{\s\s\s  http://www.maths.warwick.ac.uk/agt/AGTVol\thevolumenumber/agt-\thevolumenumber-\thepapernumber.abs.html}
\immediate\write\gtoutfile{\noexpand\\}
\immediate\write\gtoutfile{}
\ifx\theasciiabstract\relax
\immediate\write\gtoutfile{\theabstract}\else
\immediate\write\gtoutfile{\theasciiabstract}\fi
\immediate\write\gtoutfile{}
\immediate\write\gtoutfile{\noexpand\\}
\immediate\write\gtoutfile{}
\immediate\closeout\gtoutfile}}  

\def\maketitlepage{\makeagttitle\makeheadfile}
\let\makeshorttitle\maketitlepage
\let\maketitle\maketitlepage

\lognumber{42}
\volumenumber{4}
\volumeyear{2004}
\papernumber{42}
\pagenumbers{961}{1011}
\received{9 November 2003} 
\revised{20 May 2004}
\accepted{11 June 2004}
\published{31 October 2004}

\usepackage{amssymb,amsmath,amscd,graphicx}
\begin{document}

\title{A class of tight contact structures on $\Sigma_2 \times I$}
\asciititle{A class of tight contact structures on Sigma_2 x I}

\authors{Tanya Cofer}
\address{Department of Mathematics, Northeastern Illinois University\\5500 
North St Louis Avenue, Chicago, IL 60625-4699, USA}

\email{T-Cofer@neiu.edu} \urladdr{http://www.neiu.edu/~tcofer/}

\begin{abstract} 
We employ cut and paste contact topological techniques to classify
some tight contact structures on the closed, oriented genus-2
surface times the interval.  A boundary condition is specified so
that the Euler class of the of the contact structure vanishes when
evaluated on each boundary component.  We prove that there exists a
unique, non-product tight contact structure in this case.
\end{abstract}



\maketitle  

\section{Introduction}
If $M$ is a compact, oriented 3-manifold with boundary, a (positive)
\emph{contact structure} on $M$ is a completely non-integrable
2-plane distribution $\xi$ given as the kernel of a non-degenerate
1-form $\alpha$ such that $\alpha \wedge d\alpha
> 0$ at every point of $M$.  We say $\xi$ is \emph{tight} if
there is no embedded disk $D \subset M^3$ with the property that
$\xi$ is everywhere tangent to $D$ along $\partial D$.  Such a $D$
is called an \emph{overtwisted disk} and contact structures
containing such disks are called \emph{overtwisted contact
structures}.

The field of contact topology has changed profoundly and developed
rapidly during the last decades of the twentieth century. In the
1970's, Lutz and Martinet \cite{Mart} showed that every closed,
orientable three-manifold admits a contact structure. By the
1980's and early 1990's, results of Bennequin \cite{Benn} and
Eliashberg \cite{Eliash} were indicating the existence of a
qualitative difference between the classes of tight and
overtwisted contact structures. It was Eliashberg who made it
clear that the topologically interesting case to study is tight
contact structures \cite{EliashOTCS}.  He did this by showing
that, up to isotopy, overtwisted contact structures are in
one-to-one correspondence with homotopy classes of 2-plane fields
on $M$. Soon after, he gave us the classifications of tight
contact structures on $B^3$ (a foundational result for the
classification of contact structures on three manifolds), $S^3$,
$S^2 \times S^1$ and $\mathbb{R}^3$ \cite{Eliash}.  A rush of
further classification studies ensued, including the
classification of tight contact structures on the 3-torus
\cite{KandaTor}, lens spaces \cite{EtnyreLensSp,GirLensSp,TCSI},
solid tori, $T^2 \times I$ \cite{TCSI,ML}, torus bundles over
circles \cite{GirLensSp,TCSII}, and circle bundles over closed
surfaces \cite{GirCircBund,TCSII}. Etnyre and Honda
\cite{Nonexist} made a significant contribution by exhibiting a
manifold that carries no tight contact structure whatsoever.

This critical mass of understanding has recently yielded a coarse
classification principle for tight contact structures on
3-manifolds. Collectively, work by Colin, Giroux, Honda, Kazez and
Mati\'{c} \cite{atoroidalC, atoroidal2, convexdecomp} indicates
that if $M$ is closed, oriented and irreducible, then $M$ supports
finitely many isotopy classes of tight contact structures if and
only if $M$ is atoroidal.  In order to gain further understanding
of the tight contact structures on atoroidal manifolds with
infinite fundamental group, Honda Kazez and Mati\'{c} studied
hyperbolic 3-manifolds that fiber over the circle \cite{FHM}.  In
their case $M = \Sigma_g \times I$, where $\Sigma_g$ is the
genus-g surface, $I = [0,1]$ and $\xi$ is ``extremal'' with
respect to the Bennequin inequality. That is, $|e(\xi)[\Sigma_g
\times \{t\}]| = 2g-2$, $t \in [0,1]$.

In the present paper, $\Sigma_2$ is a closed, oriented genus-2
surface and $M = \Sigma_2 \times I$. We classify tight contact
structures on $(M,\Gamma_{\partial M},\mathcal{F})$ where
$\Gamma_{\partial M}$ is specified to be a single, nontrivial
separating curve on each boundary component and $\mathcal{F}$ is a
foliation which is \emph{adapted} to this dividing set.  This is a
basic case of a boundary condition satisfying $e(\xi)[\Sigma_2
\times \{0\}] = e(\xi)[\Sigma_2 \times \{1\}] = 0$.

Our classification exploits cut-and-paste methods developed by
Honda, Kazez and Mati\'{c} \cite{convexdecomp,FHM} for
constructing tight contact structures on Haken manifolds.  The
first step is to perform a Haken decomposition of $M$ in the
contact category (called a \emph{convex decomposition}). By
keeping track of certain curves (\emph{dividing curves}) on
$\partial M$ and all cutting surfaces, we specify the contact
structure in the complement of a union of 3-balls.  Under certain
conditions, we may then extend this contact structure to the
interior of each ball so that the resulting structure is tight on
the cut-open manifold. Since this type of decomposition can
generally be done in a number of ways, we apply gluing theorems to
determine which of these decompositions are associated to distinct
tight contact structures on $M$.

Our first splitting of $M$ is along a convex annulus $A$ that
separates $M$ into the disjoint union of two genus-2 handlebodies.
Although there are an infinite number of possible dividing sets on
$A$, a series of reduction arguments allow us to consider only four.
We then apply the gluing/classification theorem \cite{gluing} to a
convex decomposition of each handlebody with each of the four
dividing sets on $A$ in turn.  These convex decompositions also
allow us to locate bypasses along $A$, establishing equivalence
among some of the contact structures supported on $M \backslash A$.
This process yields an upper bound of two for the number of tight
contact structures on $M$.

At this point, it is necessary to decide whether or not the two
tight contact structures on the split-open manifold $M \backslash A$
are associated with distinct tight contact structures on $M$. Since
we must use a state-transition argument and the
gluing/classification theorem in the first stage of the
classification, we are unable to establish universal tightness for
one of the two structures. This presents an obstacle to applying the
gluing theorem directly. Thus, stage two of the classification is to
adapt the gluing theorem \cite{colgl,convexdecomp} and conclude the
existence of a unique, non-product tight contact structure on $M$.
This adaptation involves exploiting what we know about which
bypasses exist along $A$ and constructing an infinite class of
covers of $M \backslash A$ which we prove to be tight.

This process culminates in Lemma \ref{Decomp2T21POTTGT} and Theorem
\ref{Decomp2GLUING} which constitute the main focus of this paper.
They are summarized in the following statement:

\newtheorem{thm}{Theorem}[subsection]
\begin{thm}[Main Theorem]  There exists exactly one non-product tight
contact structure on  $(M,\Gamma_{\partial M},\mathcal{F})$ where
$\Gamma_{\partial M}$ is specified to be a single, nontrivial
separating curve on each boundary component and $\mathcal{F}$ is a
foliation which is \emph{adapted} to this dividing set.
\end{thm}

Once this theorem is established, we then demonstrate a special
property of the unique, non-product tight contact structure $\xi$ on
$M$. If $(M,\xi)$ is contained in some $(M^\prime,\xi^\prime)$, then
for any convex surface with boundary $S \subset M^\prime$ such that
$\partial S \subset \partial M$ and $\# (\partial S \cap
\Gamma_{\partial M}) = 2$, the dividing set on $S$ cannot contain
any boundary-parallel dividing arcs. This is because complementary
bypasses always exist inside of $(M,\xi)$.

In the following sections, we describe the background results, tools
and methods necessary for our classification.

\section{Background and tools}

\subsection{Convex surfaces}

We say that a curve $\gamma$ inside a contact manifold $(M,\xi)$ is
\emph{Legendrian} if it is everywhere tangent to $\xi$. Consider a
properly embedded surface $S \subset (M,\xi)$. Generically, the
intersection $\xi_p \cap T_p S$ at a point $p \in S$ is a vector
$X(p)$.  Integrating the vector field $X$ on $S$ gives us the a
singular foliation called the \emph{characteristic foliation}
$\xi_S$.  The leaves of the characteristic foliation are Legendrian
by definition.  A surface $S \subset M$ is called \emph{convex} if
there exists a vector field $v$ transverse to $S$ whose flow
preserves the contact structure $\xi$. Such a vector field is called
a \emph{contact vector field}.   Given a convex surface S, we define
the \emph{dividing set} $\Gamma_S = \{x \in S | v(x) \in \xi(x)\}$.
Generically, this is a collection of pairwise disjoint, smooth
closed curves (dividing curves) on a closed surface $S$ or a
collection of curves and arcs if $S$ has boundary (where the
dividing arcs begin and end on $\partial S$) \cite{convexity}. The
curves and arcs of the dividing set are transverse to the
characteristic foliation and, up to isotopy, this collection is
independent of the choice of contact vector field $v$.  Moreover,
the union of dividing curves divides $S$ into positive and negative
regions $R_{\pm}$.  $R_+ (R_-) \subset
\partial M$ is the set of points where the orientation of $\xi$
agrees (disagrees) with the orientation of $S$.

Giroux proved that a properly embedded closed surface in a contact
manifold can be $C^\infty$-perturbed into a convex surface
\cite{convexity}.  We will refer to this as the \emph{perturbation
lemma}.  If we want to keep track of the contact structure in a
neighborhood of a convex surface, we could take note of the
characteristic foliation.  However, the characteristic foliation
is very sensitive to small perturbations of the surface.
\emph{Giroux Flexibility} \cite{convexity} highlights the
usefulness of the dividing set $\Gamma_S$ by showing us that
$\Gamma_S$ captures all of the important contact topological
information in a neighborhood of $S$.  Therefore, we can keep
track of the dividing set instead of the exact characteristic
foliation.

Given a singular foliation $\mathcal{F}$ on a convex surface $S$, a
disjoint union of properly embedded curves $\Gamma$ is said to
\emph{divide} $\mathcal{F}$ if there exists some $I$-invariant
contact structure $\xi$ on $S \times I$ such that $\mathcal{F} =
\xi|_{S \times \{0\}}$ and $\Gamma$ is the dividing set for $S
\times \{0\}$.

\begin{thm}[Giroux Flexibility]
Let $S \subset (M^3,\xi)$ be a convex surface in a contact manifold
which is closed or compact with Legendrian boundary. Suppose $S$ has
characteristic foliation $\xi|_S$, contact vector field $v$ and
dividing set $\Gamma_S$.  If $\mathcal{F}$ is another singular
foliation on $S$ which is \emph{divided} by $\Gamma$, then there is
an isotopy $\phi_s$, $s$ in $[0,1]$ of $S$ such that

\begin{enumerate}
\item $\phi_0(S) = S$, \item $\xi|_{\phi_1(S)} = \mathcal{F}$, \item
$\phi$ fixes $\Gamma$, \item $\phi_s(S)$ is transverse to $v$ for all
$s$.
\end{enumerate}

\end{thm}

The \emph{Legendrian realization principle} specifies the conditions
under which a collection of curves on a convex surface $S$ may be
realized as a collection of Legendrian curves by perturbing $S$ to
change the characteristic foliation on $S$ while keeping $\Gamma_S$
fixed.  In general this is not a limiting condition. The result is
achieved by isotoping the convex surface $S$ through surfaces that
are convex with respect to the contact vector field $v$ for $S$ so
that the collection of curves on the isotoped surface are made
Legendrian.

Let $C$ be a collection of closed curves and arcs on a convex
surface $S$ with Legendrian boundary.  We call $C$
\emph{nonisolating} if:
\begin{enumerate}
    \item $C$ is transverse to $\Gamma_S$.
    \item Every arc of $C$ begins and ends on $\Gamma_S$.
    \item The elements of $C$ are pairwise disjoint.
    \item If we cut $S$ along $C$, each component intersects the dividing set $\Gamma_S$.
\end{enumerate}
An isotopy $\phi_s$, $s \in [0,1]$ of a convex surface $S$ with
contact vector field $v$ is called \emph{admissible} if $\phi_s(S)$
is transverse to $v$ for all $s$.

\begin{thm} [Legendrian Realization]
If $C$ is a nonisolating collection of disjoint, properly embedded
closed curves and arcs on a convex surface  $S$ with Legendrian
boundary, there is an admissible isotopy $\phi_s$, $s \in [0,1]$ so
that:
\begin{enumerate}
    \item $\phi_0 = id$,
    \item each surface $\phi_s(S)$ is convex,
    \item $\phi_1(\Gamma_S)$ = $\Gamma_{\phi_1(S)}$,
    \item $\phi_1(C)$ is Legendrian.
\end{enumerate}
\label{legreal}
\end{thm}

A useful corollary to Legendrian realization was formulated by
Kanda \cite{KandaLR}:
\newtheorem{corollary}{Corollary}[section]
\begin{corollary}
Suppose a closed curve $C$ on a convex surface $S$
\begin{enumerate}
    \item is transverse to $\Gamma_S$,
    \item nontrivially intersects $\Gamma_S$.
\end{enumerate}
Then $C$ can be realized as a Legendrian curve. \label{legrealK}
\end{corollary}

\begin{figure}[ht!]
\begin{center}
\includegraphics[bb=0 0 108 104]{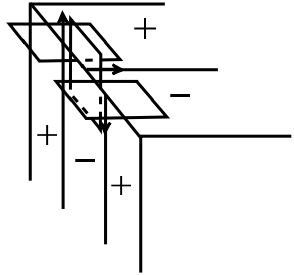}
\end{center}
\caption{The dividing set on a cutting surface}
\label{TwistingPlanes}
\end{figure}

Suppose $\gamma \subset S$ is Legendrian and $S \subset M$ is a
properly embedded convex surface. We define the
\emph{Thurston-Bennequin number} $tb(\gamma, Fr_S)$ of $\gamma$
relative to the framing, $Fr_S$, of $S$ to be the number of full
twists $\xi$ makes relative to $S$ as we traverse $\gamma$, where
left twists are defined to be negative.  It turns out that
$tb(\gamma, Fr_S)$ = $-\frac{1}{2}\#(\Gamma_S \cap \gamma)$ (see
figure \ref{TwistingPlanes}). When $\gamma$ is not a closed curve,
we will refer to the \emph{twisting} $t(\gamma, Fr_S)$ of the arc
$\gamma$ relative to $Fr_S$.

Given any Legendrian curve $\gamma = \partial S$ with non-positive
Thurston-Bennequin number, the following relative version of
Giroux's perturbation lemma, proved by Honda \cite{TCSI}, asserts
that we can always arrange for the contact planes to twist
monotonically in a left-handed manner as we traverse the curve.
That is, an annular neighborhood $A$ of a Legendrian curve
$\gamma$ with $tb(\gamma) = -n$ inside a contact manifold
($M$,$\xi$) is locally isomorphic to $\{x^2 + y^2 \leq  \epsilon
\} \subset (\mathbb{R}^2 \times \mathbb{R} / \mathbb{Z}, (x,y,z))$
with contact 1-form $\alpha = \sin(2 \pi n z)dx + cos(2 \pi n
z)dy$, $n \in \mathbb{Z}^+$ \cite{TCSI,gluing}.  This is called
\emph{standard form}.  Once this is achieved, it is then possible
to perturb the surface $S$ so as to make it convex.

\begin{thm} [Relative Perturbation] Let $S \subset M$ be a compact, oriented, properly embedded surface with Legendrian boundary such that $t(\gamma,Fr_S) \leq 0$ for all components $\gamma$ of $\partial S$.  There exists a $C^0$- small perturbation near the boundary which fixes $\partial S$ and puts an annular neighborhood $A$ of $\partial S$ into standard form.  Then, there is a further $C^\infty$-small perturbation (of the perturbed surface, fixing $A$) which makes $S$ convex.  Moreover, if $v$ is a contact vector field on a neighborhood of $A$ and transverse to $A$, then $v$ can be extended to a vector field transverse to all of $S$.
\label{convexthm}
\end{thm}

\subsection{The method: decomposing $(M,\xi)$}

We will be analyzing $(M^3,\xi)$ by cutting it along surfaces into
smaller pieces and then analyzing the pieces.  A familiar way to do
this is to cut down $M$ along a sequence of incompressible surfaces
$\{S_i\}$ into a union of balls.  This is known as a \emph{Haken
decomposition}:  $$M = M_0 \stackrel{S_1}{\leadsto} M_1
\stackrel{S_2}{\leadsto} \cdots \stackrel{S_n}{\leadsto} M_n =
\amalg B^3$$ In order to do define and exploit an analogous
procedure in the contact category, we will need to make use of a
couple of important results, including the following theorem, due to
Eliashberg, which is a foundational result for the cut-and-paste
methods described here:

\begin{thm}[Eliashberg's Uniqueness Theorem]
If $\xi$ is a contact structure in a neighborhood of $\partial B^3$
that makes $\partial B^3$ convex and the dividing set on $\partial
B^3$ consists of a single closed curve, then there is a unique
extension of $\xi$ to a tight contact structure on $B^3$ (up to
isotopy that fixes the boundary).
\end{thm}

Also,  Giroux proved that we can determine when convex surfaces have
tight neighborhoods by looking at their dividing sets
\cite{convexity}.

\begin{thm}[Giroux's Criterion]
Suppose $S \neq S^2$ is a convex surface in a contact manifold
$(M,\xi)$.  There exists a tight neighborhood for $S$ if and only if
$\Gamma_S$ contains no homotopically trivial closed curves.  If $S =
S^2$, then $S$ has a tight neighborhood if and only if
$\#\Gamma_{S^2} = 1$.
\end{thm}

By \emph{Giroux's criterion}, we know that if a homotopically
trivial dividing curve appears on a convex surface inside a contact
manifold, then our contact structure is overtwisted.  If our surface
happens to be $S^2$, then more than one dividing curve indicates an
overtwisted contact structure.

Let $M$ be an irreducible contact 3-manifold with boundary.  Then
$M$ admits a Haken decomposition along incompressible surfaces
$\{S_i\}$ \cite{haken}.  To do this in the contact category, we may
take one of two approaches.  In one approach, we assume $M$ carries
a contact structure $\xi$ which makes $\partial M$ convex. In this
case we may, at each stage of the decomposition, (1) use the
relative perturbation lemma to perturb $S_i$ into a convex surface
with Legendrian boundary, and (2) cut along $S_i$ and \emph{round
edges} (see the following discussion) so that the cut-open manifold
is a smooth contact manifold with convex boundary \cite{TCSTF}.

Alternatively, we can do this decomposition abstractly along
surfaces with divides with the goal of \emph{constructing} a tight
contact structure on $M$ (which may or may not exist).  In this
case, we begin with a \emph{convex structure}
($M,\Gamma,R_+(\Gamma),$ $R_-(\Gamma)$) \cite{TCSTF} where $\Gamma$,
$R_+(\Gamma)$ and $R_-(\Gamma)$ satisfy all the properties of a
dividing set and positive and negative regions on $\partial M$ were
$\partial M$ to be convex with respect to some contact structure
$\xi$ on $M$. Then we can cut open $M$ along surfaces with divides
$S_i$ that we hope will be convex with respect to the contact
structure that we are attempting to construct.  Note that we will
make no distinction between these ``\emph{pre-convex}'' surfaces and
actual ones since the contact structure in a (contact) product
neighborhood of the surface is determined up to isotopy by its
dividing set. If we are able to decompose $M$ in this manner and we
end up with a union of $(B^3,S^1)$'s, then we may use Giroux
flexibility and Eliashberg's uniqueness theorem to conclude the
existence of a contact structure $\xi$ on $M$ which is compatible
with the decomposition, and which is tight on the cut-open manifold.
If there is a tight contact structure on $M$ which makes $\partial
M$ convex, it is always possible to perform a convex decomposition
of $(M,\Gamma)$ which is compatible with $\xi$ \cite{TCSTF}.  Hence,
tight contact structures are determined by their convex
decompositions.

Before we define this decomposition properly, we will explain what
happens to the dividing sets on convex surfaces in a cut-open
contact manifold when we smooth out the edges. We make use of the
\emph{edge rounding} lemma, a result that tells us what happens to
the contact structure in a neighborhood of a Legendrian curve of
intersection when we smooth away corners. A proof of the edge
rounding lemma can be found in \cite{TCSI}. A useful and brief
discussion of the result can be found in \cite{gluing}. The
following is a summary.

\begin{figure}[ht!]
\begin{center}
\includegraphics[bb=0 0 230 104]{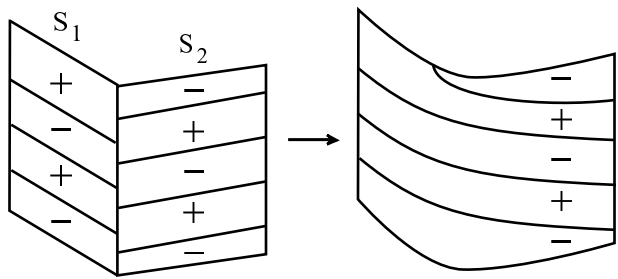}
\end{center}
\caption{Edge Rounding} \label{EdgeRounding}
\end{figure}

If we consider two compact, convex surfaces $S_1$ and $S_2$ with
Legendrian boundary that intersect along a common Legendrian
boundary curve, then a neighborhood of the common boundary is
locally isomorphic to $\{x^2 + y^2 \leq  \epsilon \} \subset
(\mathbb{R}^2 \times \mathbb{R} / \mathbb{Z}, (x,y,z))$ with contact
1-form $\alpha = \sin(2 \pi n z)dx + cos(2 \pi n z)dy$, $n \in
\mathbb{Z}^+$ (standard form).  We let $A_i \subset S_i$, $i \in
\{1,2\}$ be an annular collar of the common boundary curve. It is
possible to choose this local model so that $A_1 = \{x=0, 0 \leq y
\leq \epsilon \}$ and $A_2 = \{y=0, 0 \leq x \leq \epsilon \}$.
Then the two surfaces are joined along $x = y = 0$ and rounding this
common edge results in the joining of the dividing curve $z =
\frac{k}{2n}$ on $S_1$ to $z = \frac{k}{2n} - \frac{1}{4n}$ on $S_2$
for $k \in \{0,...,2n-1\}$.  See figure \ref{EdgeRounding}.

We are now ready to define our decomposition technique.  We will
use an abstract formulation for a convex decomposition in our
classification.  To do this, we will apply theorems such as
Legendrian realization, relative perturbation, and edge rounding
abstractly, taking care that our abstract divides satisfy the
hypotheses of these theorems and the effect of applying of these
theorems to the surfaces and divides follows as it would in the
presence of an appropriate contact structure on $M$.

Whenever it is possible to find such a decomposition into a union of
3-balls $B^3$ with $\# \Gamma_{\partial B^3} = 1$ , we say that
$$(M,\Gamma) = (M_0,\Gamma_0) \stackrel{(S_1,\sigma_1)}\leadsto (M_1,\Gamma_1) \stackrel{(S_12,\sigma_2)}\leadsto \cdots \stackrel{(S_n,\sigma_n)}\leadsto (M_n,\Gamma_n) = \amalg (B^3,S^1)$$
defines a \emph{convex decomposition} of $(M,\Gamma)$ where the
$\sigma_i$ are a set of divides on the convex surface $S_i$ and
$M_{i} = M_{i-1} \backslash S_i$ inherits $\Gamma_{i}$ from
$\Gamma_{i-1}$, $\sigma_{i}$ and (abstract) edge rounding.

Now, we would like to know how $\Gamma_{S_i}$ changes if we cut
along a different, but isotopic, surface $S_i^\prime$ with the same
boundary inside a contact manifold $(M,\xi)$. We may then return to
our notion of a convex decomposition and apply the results
abstractly to our cutting surfaces with divides. Honda and Giroux
\cite{GirLensSp,TCSI} both studied the changes in the characteristic
foliation on a convex surface under isotopy. Honda singled out the
minimal, non-trivial isotopy of a cutting surface which he calls a
\emph{bypass}.

A \emph{bypass} (figure \ref{bypassdisk}) for a convex surface $S
\subset M$ (closed or compact with Legendrian boundary) is an
oriented, embedded half-disk $D$ with Legendrian boundary such
that $\partial D$ = $\gamma_1 \cup \gamma_2$ where $\gamma_1$ and
$\gamma_2$ are two arcs that intersect at their endpoints.  $D$
intersects $S$ transversely along $\gamma_1$ and $D$ (with either
its given orientation or the opposite one) has positive elliptic
tangencies along $\partial \gamma_1$, a single negative elliptic
tangency along the interior of $\gamma_1$, and only positive
tangencies along $\gamma_2$, alternating between elliptic or
hyperbolic. Moreover, $\gamma_1$ intersects $\Gamma_S$ in exactly
the three elliptic points for $\gamma_1$ \cite{TCSI}.

\begin{figure}[ht!]
\begin{center}
\includegraphics[bb=0 0 210 111]{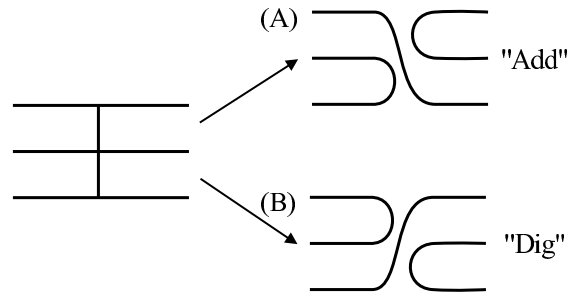}
\end{center}
\caption{Abstract Bypass Moves} \label{absbypassmove}
\end{figure}

Isotoping a cutting surface past a bypass disk (figure
\ref{bypassdisk}) produces a change in the dividing set as shown
in figure \ref{absbypassmove}.  Figure \ref{absbypassmove}(A)
shows how the dividing set changes if we attach a bypass above the
surface. Figure \ref{absbypassmove}(B) shows the change in the
dividing set if we dig out a bypass below the surface. Formally,
we have:

\begin{figure}[ht!]
\begin{center}
\includegraphics[bb=0 0 140 101]{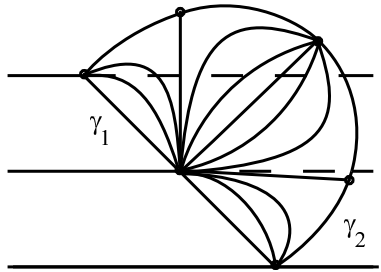}
\end{center}
\caption{Bypass Disk} \label{bypassdisk}
\end{figure}

\begin{thm}[Bypass Attachment]
Suppose $D$ is a bypass for a convex $S \subset M$.  There is a
neighborhood of  $S \cup D$ in $M$ which is diffeomorphic to $S
\times [0,1]$ with $S_i$ = $S \times \{i\}$, $i \in \{0,1\}$ convex,
$S \times [0,\epsilon]$ is $I$-invariant, $S = S \times
\{\epsilon\}$ and $\Gamma_{S_1}$ is obtained from $\Gamma_{S_0}$ as
in figure \ref{absbypassmove} (A).
\end{thm}

A useful result concerning bypass attachment is the \emph{bypass
sliding lemma} \cite{TCSII,FHM}.  It allows us some flexibility with
the choice of the Legendrian arc of attachment.

Let $C$ be a curve on a convex surface $S$ and let $M =
min\{\#(C^\prime \cap \Gamma_S)|$ $C^\prime$ is isotopic to $C$ on
$S\}$. We say that $C$ is \emph{efficient} with respect to
$\Gamma_S$ if $M \neq 0$ and the geometric intersection number $\#(C
\cap \Gamma_S) = M$, or, if $M = 0$, then $C$ is \emph{efficient}
with respect to $\Gamma_S$ if $\#(C \cap \Gamma_S) = 2$.

\begin{thm}[Bypass Sliding Lemma]
Let $R$ be an embedded rectangle with consecutive sides $a,b,c$ and
$d$ on a convex surface $S$ so that $a$ is the arc of attachment of
a bypass and $b,d \subset \Gamma_S$.  If $c$ is a Legendrian arc
that is \emph{efficient} (rel endpoints) with respect to $\Gamma_S$,
then there is a bypass for which $c$ is the arc of attachment.
\end{thm}

\subsection{Gluing}

When we do a convex decomposition of $(M^3,\Gamma)$ and end up with
a union of balls, each with a single dividing curve, we know we have
a contact structure on $M$ which is tight on this cut-open manifold.
If we glue our manifold back together along our cutting surfaces,
the contact structure may not stay tight.  It may be that an
overtwisted disk $D \subset M$ intersected one or more of the
cutting surfaces.  There are two gluing theorems that tell us when
we can expect the property of tightness to survive the regluing
process.

One gluing theorem, which we will discuss in more detail later, is
due to Colin \cite{colgl}, Honda, Kazez and Mati\'{c}
\cite{convexdecomp}. This theorem allows us to conclude tightness of
the glued-up manifold when the dividing sets on the cutting surfaces
are \emph{boundary-parallel}.  That is, the dividing set consists of
arcs which cut off disjoint half-disks along the boundary. We will
see that we must adapt this theorem in order to complete the
classification in the case of our work here.

\begin{thm}[Gluing]
Consider an irreducible contact manifold $(M,\xi)$ with nonempty
convex boundary and $S \subset M$ a properly embedded, compact,
convex surface with nonempty Legendrian boundary such that:
\begin{enumerate}
    \item $S$ is incompressible in $M$,
    \item $t(\delta,Fr_S) < 0$ for each connected component $\delta$ of $\partial S$ (i.e.\ each $\delta$ intersects $\Gamma_{\partial M}$ nontrivially), and
    \item $\Gamma_S$ is boundary-parallel.
\end{enumerate}
If we have a decomposition of $(M,\xi)$ along $S$, and $\xi$ is
universally tight on $M \backslash S$, then $\xi$ is universally
tight on $M$. \label{GLUING}
\end{thm}

The other gluing theorem is due to Honda \cite{gluing}. Honda
discovered that, if we take a convex decomposition of an
overtwisted contact structure on $M$ and look at all possible
non-trivial isotopies (bypasses) of the cutting surfaces $S_i$, we
will eventually come upon a decomposition that does not cut
through the overtwisted disk.

In order to partition the set of contact structures on $M$ into
isotopy classes of tight and overtwisted structures, we appeal to a
mathematical algorithm known as the gluing/classification theorem.
We describe this algorithm in the case where $M$ is a handlebody and
$\xi$ is a contact structure on $M$ so that $\partial M$ is convex,
since it is this case that is most relevant for our work here.  One
can find a statement of the general case in Honda \cite{gluing}.

Let $M = H_g$ be a handlebody of genus $g$.  Prescribe
$\Gamma_{\partial H_g}$ and a compatible characteristic foliation.
Let $\{D_i\}$, $i \in \{ 1...g \}$ be a collection of disjoint
compressing disks with $tb(\partial D_i, Fr_{D_i}) < 0$ yielding the
convex splitting $M = H_g \leadsto H_g \backslash (D_1 \cup D_2 \cup
\cdots \cup D_g) = B^3 = M^\prime$. Also, let $C =
(\Gamma_1,\Gamma_2, \cdots ,\Gamma_g)$ represent a configuration of
dividing sets on these disks, which we will call a \emph{state}.  We
now need to decide if a given state corresponds to a tight contact
structure on the original manifold $M$.  We call a state $C$ {\em
potentially allowable} (i.e.\ not obviously overtwisted) if
($M$,$\Gamma_{\partial M}$), cut along $D_1 \cup D_2 \cup \cdots
\cup D_g$ with configuration $C$, gives the boundary of a tight
contact structure on $B^3$ (after edge rounding). That is,
$\Gamma_{S^2}$ consists of a single closed curve.

Honda introduced the idea of the \emph{state transition} and defined
an equivalence relation on the set of configurations $\mathcal C_M$
on a manifold $M$ \cite{gluing}.  This equivalence relation
partitions  $\mathcal C_M$ into equivalence classes so that each
equivalence class represents a distinct contact structure on $M$,
either tight or overtwisted.

We say that a {\em state transition} $C \stackrel{st}{\leadsto}
C^\prime$ exists from a state $C$ to another state $C^\prime$ on $M$
provided:

\begin{enumerate}
\item $C$ is \emph{potentially allowable}. \item There is a
nontrivial abstract bypass dig (see figure \ref{absbypassmove}) from
one copy of some $D_i$ in the cut-open $M^\prime$ and a
corresponding add along the other copy of $D_i$ that produces
$C^\prime$ from $C$ (a trivial bypass does not alter the dividing
curve configuration). \item  Performing only the dig from item (2)
(above) does not change the number of dividing curves in
$\Gamma_{\partial B^3}$ (i.e.\ the bypass disk exists inside of $M
\backslash C$).
\end{enumerate}

Define $<\sim>$ on $\mathcal C_M$ as the equivalence relation
generated by $\sim$, where $C \sim C^\prime$ if either $C
\stackrel{st}{\leadsto} C^\prime$ or $C^\prime
\stackrel{st}{\leadsto} C$. Then, an equivalence class under
$<\sim>$ represents a tight contact structure provided every state
in that class is potentially allowable (i.e.\ no state in the
equivalence class is obviously overtwisted).

\begin{thm}[Gluing/Classification]
The tight contact structures on ($H_g$, $\Gamma_{\partial H_g}$)
are in one-to-one correspondence with the equivalence classes
under $<\sim>$ containing only potentially allowable states.
\label{gluingclassification}
\end{thm}

\subsection{State transition arguments}

The gluing/classification theorem for the case of $M = H_g$, as
stated above, is reasonably practical.  This is due to the fact
that, although the number may be large, there are at least a finite
number of possible states to check. In general, straightforward
applications of gluing/classification theorem are impractical due to
the possibility of an infinite number of dividing curve
configurations (\emph{states}) on general cutting surfaces. However,
state transition arguments may still be used in a more general
setting to show two contact structures are equivalent. Consider the
case of our work here: $M = \Sigma_2 \times I$ with fixed
$\Gamma_{\partial M}$. Our classification begins by cutting $M$
along a convex annulus that separates $M$ into two genus-2
handlebodies.  The infinite number of possible dividing curve
configurations on the annulus precludes the use of
 the gluing/classification theorem on $M$.  However, given
two dividing curve configurations on the annulus, say $\Gamma_{A}$ =
$A_1$ and $\Gamma_A$ = $A_2$, we may conclude that (1) both
configurations give the boundary of a tight contact structure on the
cut-open manifold $M \backslash A$, and (2) contact structures on
$(M \backslash A, A_1)$ and $(M \backslash A, A_2)$ come from the
same contact structure on $M$, just cut along different, but
isotopic annuli with the same boundary.  This is achieved as
follows:
\begin{enumerate}
    \item Consider contact structures on $(M \backslash A,
A_1)$ and $(M \backslash A, A_2)$ specified by a convex
decomposition of each, and suppose the gluing/classification theorem
shows that the contact structure on $(M \backslash A, A_1)$ is
tight.
    \item Consider the same cutting disks $D_i$ and $D^{'}_i$, $i{\in}\{1,2\}$ for both $(M \backslash
A, A_1)$ and $(M \backslash A, A_2)$ with the contact structures
specified above given by a choice of dividing set on each cutting
disk. If there exists a bypass $B$ along $A_1$ (we must check this)
so that isotoping past $B$ transforms $A_1$ into $A_2$ and
simultaneously changes $\Gamma_{D_i}$ into $\Gamma_{D_i^{'}}$ for
each $i$, then we say that there is a \emph{state transition} taking
$(A_1,\Gamma_{D_1},\Gamma_{D_2})$ to
$(A_2,\Gamma_{D_1^{'}},\Gamma_{D_2^{'}})$.  Any two contact
structures related by a sequence of state transitions represent the
same contact structure on $M$ just cut along different, but isotopic
annuli with the same boundary.  Note that this process establishes
equivalence of contact structures on $M$, but not tightness.  This
issue will be dealt with separately.
\end{enumerate}

\subsection{Some special bypasses}

\begin{figure}[ht!]
\begin{center}
\includegraphics[bb=0 0 189 40]{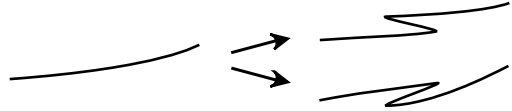}
\end{center}
\caption{Decreasing $tb(\gamma)$} \label{decreasingtb2}
\end{figure}

Establishing the existence of bypasses requires work.  In general,
the bypasses along a convex surface $S$ inside a contact manifold
$(M,\xi)$ may be ``long'' or ``deep'' (i.e.\ they exist outside of
an $I$-invariant neighborhood), and establishing existence
requires global information about the ambient manifold $M$.  This
is further evidenced by examining that the two Legendrian arcs,
$\gamma_1$ and $\gamma_2$ (where $\gamma_1 \subset S$) that
comprise the boundary of the bypass half-disk $B$ (see figure
\ref{bypassdisk}).  When we isotope the convex surface $S$ past
$B$ to produce a new convex surface $S^\prime$, we see that
$t(\gamma_1,Fr_S) = -1$ whereas $t(\gamma_2,Fr_{S^\prime}) = 0$.
Since a small neighborhood of a point on a Legendrian curve is
isomorphic to a neighborhood of the origin in $\mathbb{R}^3$ with
the standard contact structure, the Thurston-Bennequin number can
be decreased by one of the two moves in figure
\ref{decreasingtb2}, yet Benneqin's inequality tells us that
$tb(\gamma,Fr_S)$ is bounded above by the Seifert genus of
$\gamma$.  So, although it is easy to decrease twisting number, it
not always possible to increase it.

There are two types of bypasses, however, that can be realized
locally (in an $I$-invariant neighborhood of a convex surface):
\emph{trivial} bypasses and \emph{folding} bypasses.
\emph{Trivial} bypasses are those that do not change the dividing
curve configuration.  Honda argues existence and triviality of
\emph{trivial} bypasses in \cite{gluing}, lemmas 1.8 and 1.9.
\emph{Folding} bypasses are the result of certain isotopies of a
convex cutting surface inside an $I$-invariant neighborhood of the
surface.

One way to establish existence of a bypass is to cut down the
manifold $M$ into a union of 3-balls and invoke Eliahsberg's
uniqueness theorem.  Since there is a unique tight contact structure
on the 3-ball with $\# \Gamma_{\partial B^3} = 1$ (rel boundary),
only trivial bypasses exist in this case.  Thus, if we can show that
a proposed bypass along a convex surface $S \subset M$ is trivial on
$B^3$, we can conclude it's existence.

A result related to trivial bypasses is the \emph{semi-local
Thurston-Bennequin inequality}.  It is useful when we are attempting
to distinguish between the product contact structure on $M$ and
other structures.  The following formulation is translated from
Giroux \cite{GirCircBund}:

\newtheorem{proposition}{Proposition}[section]
\begin{proposition}[Semi-Local Thurston-Bennequin Inequality]  Let
$\xi$ be an $I$-invariant (product) contact structure on the product
$U = S \times I$, where $I = [-1,1]$.  Suppose $C$ is a simple
closed curve on $S = S \times \{0\}$, and $\Gamma$ is a dividing set
on $S$ which is \emph{adapted} to $\xi_S$.  Then, for all isotopies
$\phi_t$ of $U$ which take $C$ to a Legendrian curve $\phi_1(C)$,
the number of twists $\xi$ makes along $\phi_1(C)$ relative to the
tangent planes to $\phi_1(S)$ satisfies the inequality
$$t(\xi,Fr_{\phi_1(S)},\phi_1(C)) \leq -\frac{1}{2}\#(\Gamma,C)$$
where $t(\xi,Fr_{\phi_1(S)},\phi_1(C))$ measures the twisting of
$\xi$ along $\phi_1(C)$ relative to $Fr_{\phi_1(S)}$ and
$\#(\Gamma,C)$ is $min(\Gamma \cap C^\prime)$ where the minimum is
taken over all closed curves $C^\prime$ isotopic to $C$ on $S$.
Moreover, there is an isotopy which realizes equality.
\label{localtb}
\end{proposition}

Above, we mentioned that isotoping a convex surface $S$ past a
bypass $B$ where $\partial B = \gamma_1 \cup \gamma_2$ increases the
twisting number $t(\gamma_1,Fr_S) = -1$ to
$t(\gamma_2,Fr_{S^\prime}) = 0$.  So, if $C$ (containing $\gamma_1$)
is a simple closed curve satisfying $\Gamma \cap C$ =
$\#(\Gamma,C)$, then isotoping $S$ past $B$ causes the inequality to
fail.  Thus, it must be that the product structure can contain no
non-trivial bypasses.

We say that a closed curve is \emph{nonisolating} if every component
of $S \backslash \gamma$ intersects $\Gamma_S$.  A \emph{Legendrian
divide} is a Legendrian curve such that all the points of $\gamma$
are tangencies.  The ideas for the following exposition are borrowed
from Honda et al.~\cite{convexdecomp}.

\begin{figure}[ht!]
\begin{center}
\includegraphics[bb=0 0 319 194]{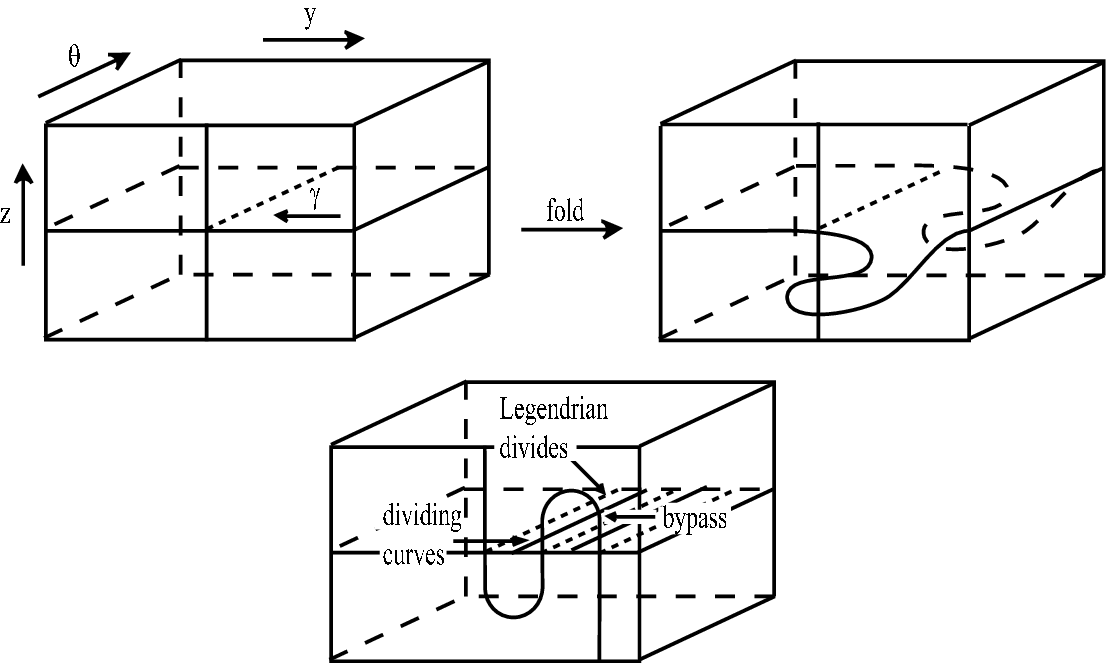}
\end{center}
\caption{Folding} \label{folding}
\end{figure}

To produce a \emph{folding} bypass, we must start with a
\emph{nonisolationg} closed curve $\gamma$ on a convex surface $S$
that does not intersect $\Gamma_S$. A strong form of
\emph{Legendrian realization} allows us to make $\gamma$ into a
\emph{Legendrian divide}.  Then, there is a local model $N = S^1
\times [-\epsilon,\epsilon] \times [-1,1]$ with coordinates
$(\theta,y,z)$ and contact form $\alpha = dz - yd\theta$ in which
the convex surface $S$ intersects the model as $S^1 \times
[-\epsilon,\epsilon] \times \{0\}$ and the Legendrian divide
$\gamma$ is the $S^1$ direction (see figure \ref{folding}).

\begin{figure}[ht!]
\begin{center}
\includegraphics[bb=0 0 296 37]{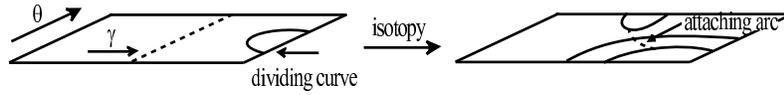}
\end{center}
\caption{Folding is isotopy plus bypass dig} \label{folding2}
\end{figure}

To \emph{fold} around the \emph{Legendrian divide} $\gamma$, we
isotope the surface $S$ into an ``S''-shape (see the figure) inside
of the model $N$. Outside of the model, $S$ and $S^\prime$ are
identical. The result is a pair of dividing curves on $S^\prime$
parallel to $\gamma$.  Note that the obvious bypass add indicated in
figure \ref{folding} ``undoes" the fold.  It turns out that we can
view the fold as an isotopy of the dividing set followed by a bypass
dig as illustrated in figure \ref{folding2}, where the bypass dig is
the one dual to the bypass add in figure \ref{folding}.

In many of the applications that follow, we will be establishing
the existence of folding bypasses on solid tori.  Note that, on
the boundary of a solid torus with $2n$ dividing curves, there are
$2n$ Legendrian divides spaced evenly between the dividing curves.
Thus, there is always a \emph{Legendrian divide} located near an
existing dividing curve.
\section{Characterizations of potentially tight contact structures on $\Sigma_2 \times I$ with a specified configuration on the boundary}

\begin{figure}[ht!]
\begin{center}
\includegraphics[bb=0 0 200 146]{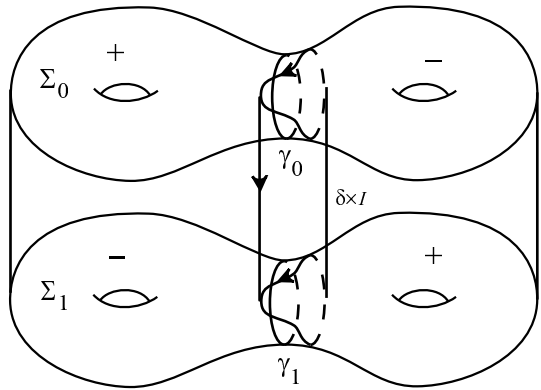}
\end{center}
\caption{M = $\Sigma_2 \times I$} \label{MisMtimesI2}
\end{figure}

We investigate the number of tight contact structures on $M
=\Sigma_2 \times I$ with a specific dividing set on the boundary
$\partial (\Sigma_2 \times I) = \Sigma_0 \cup \Sigma_1$.  The
dividing set we specify on the boundary consists of a single
separating curve $\gamma_i$ on each $\Sigma_i$, $i \in \{0,1\}$.
These curves are chosen so that $\chi((\Sigma_i)_+) =
\chi((\Sigma_i)_-)$.  Here, $(\Sigma_i)_\pm$ $i \in \{0,1\}$
represent the positive and negative regions of $\Sigma_i
\backslash \Gamma_{\Sigma_i}$.  The position of the $\gamma_i$ are
indicated in figure \ref{MisMtimesI2}.

\begin{figure}[ht!]
\begin{center}
\scalebox{0.90}{\includegraphics[bb=0 0 241 197]{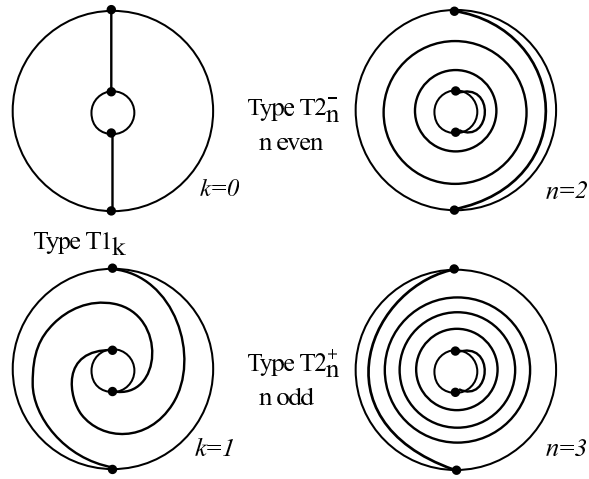}}
\end{center}
\caption{Classification of $\Gamma_A$} \label{GammaAtypes}
\end{figure}

\begin{figure}[ht!]
\begin{center}
\scalebox{0.90}{\includegraphics[bb=0 0 223 213]{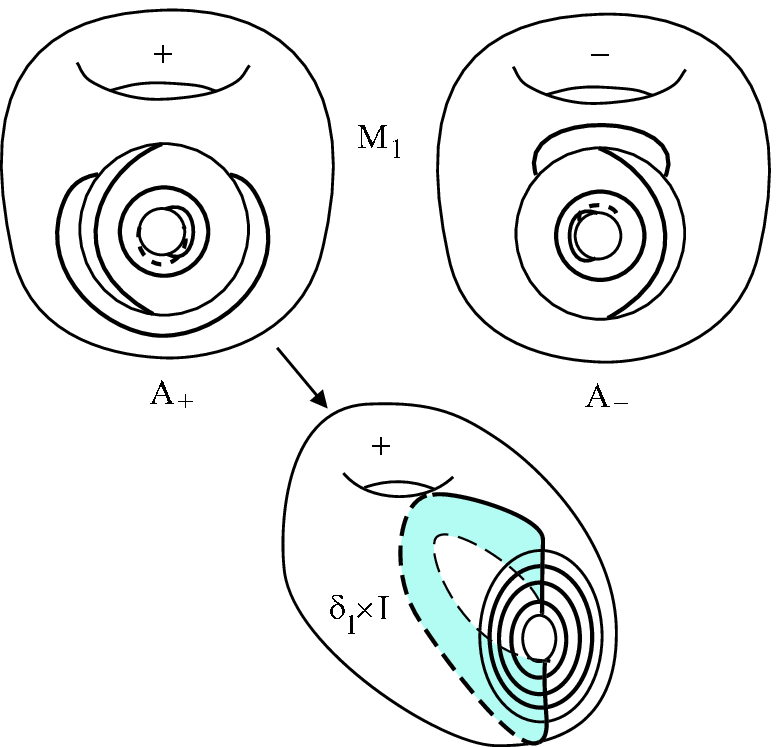}}
\end{center}
\caption{$M_1 = M \backslash A$} \label{M1isMminusA2}
\end{figure}

Our first strategy is to provide a convex decomposition of $M$ and
partition the set of contact structures into equivalence classes
using Honda's gluing/classific\-ation theorem \cite{gluing}.  One
difficulty here is that the first cutting surface ($\delta \times
I$ as indicated in figure \ref{MisMtimesI2}) of the decomposition
is an annulus $A$. Fortunately, the infinitely many possibilities
for $\Gamma_A$ fall into three natural categories (see figure
\ref{GammaAtypes}). It is immediately clear in this case that
$\Gamma_A$ of type $T2_{k}^-$ are overtwisted. For the following
reduction arguments, we will assume that there is a tight contact
structure on $M$ with $\Gamma_A$ as specified, and we prove that
we can always find a bypass along $A$ so that isotoping $A$ past
this bypass produces the desired reduction.  We can assume that
$\delta \times I$ is convex with Legendrian boundary since
$\delta$ intersects $\Gamma_{\partial M}$ nontrivially. $ M_1 = M
\backslash A = M_1^+ \cup M_1^-$ is given at the top of figure
\ref{M1isMminusA2}. We show first that all tight $T2_{2k+1}^+$,
$T1_k$, and $T1_{-k}$ for $k \geq 1$ can be reduced to $T2_1^+$,
$T1_1$, and $T1_{-1}$, respectively.

\begin{figure}[b!]
\vspace{0.5cm}
\begin{center}
\scalebox{0.85}{\includegraphics[bb=0 0 307 202]{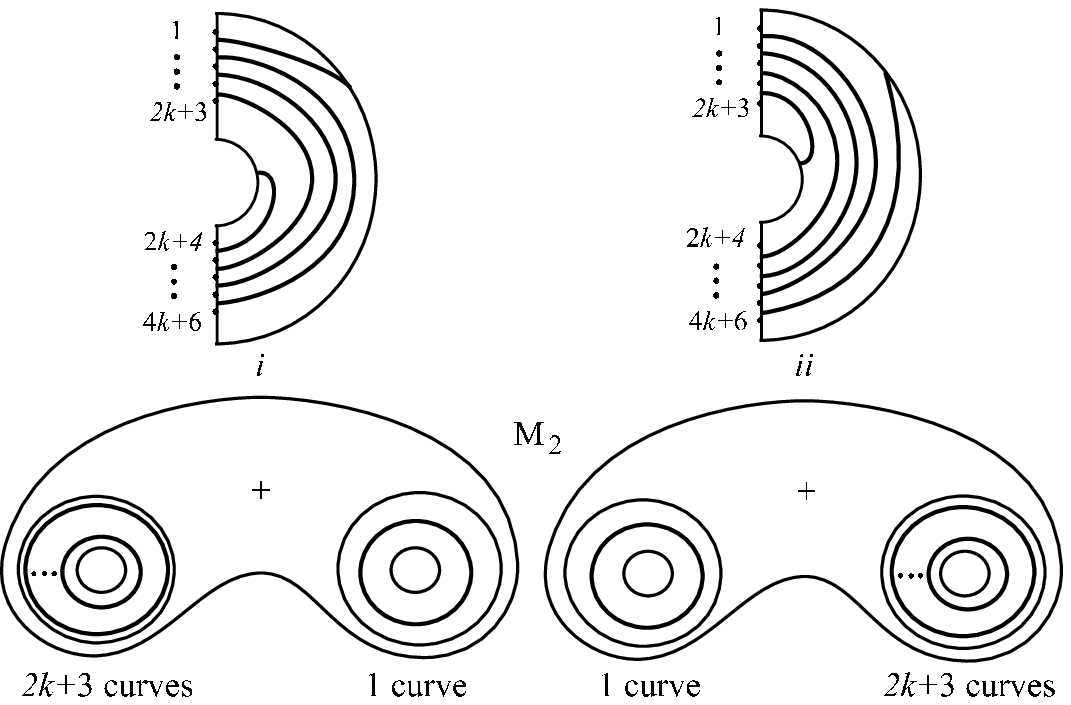}}
\end{center}
\caption{Convex decomposition for $M$ with $\Gamma_A$ =
$T2_{2k+1}^+$} \label{Decomp2T2ODDRED}
\end{figure}

\newtheorem{lem}[thm]{Lemma}
\begin{lem} $T2_{2k+1}^+$, $k \in {\bf \mathbb{Z}^+}$ can be reduced to $T2_1^+$.
\label{lemmaT2oddred2}
\end{lem}
\proof  Suppose there is a tight
contact structure on $M$ with $\Gamma_A$ = $T2_{2k+1}^+$, $k \geq
1$. Cut the $A_+$ component of $M_1$ open along the convex cutting
surface $\epsilon = \delta_1 \times I$ where $\delta_1$ is
positioned as in figure \ref{M1isMminusA2}. We see that all but
two possibilities for $\Gamma_{\epsilon_+}$ contain a dividing
curve straddling one of the positions 2 through 2\emph{k}+2 or
2\emph{k}+5 through 4\emph{k}+5 as shown in figure
\ref{Decomp2T2ODDRED}. Isotoping $A_+$ across any of these
bypasses yields a dividing set on the isotoped annulus equivalent
to $T2_{2k-1}^+$.

\begin{figure}[ht!]
\begin{center}
\includegraphics[bb=0 0 271 136]{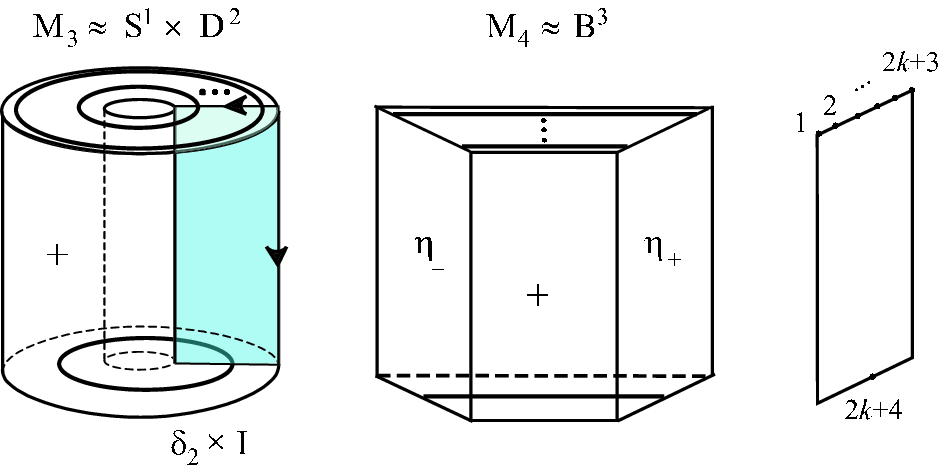}
\end{center}
\vspace{-0.3cm}
\caption{\label{T2oddred8} $\Gamma_A$ of type $T2_{2k+1}^+$}
\end{figure}

The remaining possibilities are pictured.  Both possibilities
\emph{i} and \emph{ii} give us a dividing set on $M_2 \cong S^1
\times D^2 = M_1^+ \backslash \epsilon$ consisting of 2\emph{k}+4
longitudinal curves as shown in figure \ref{Decomp2T2ODDRED}. If we
choose a convex, meridional cutting surface $\eta$, we see, as in
figure \ref{T2oddred8}, that all possibilities for $\Gamma_{\eta_+}$
contain a dividing curve straddling one of the positions 2 through
2\emph{k}+ 2 for $k \geq 1$. Isotoping $A_+$ across any of these
bypasses yields a dividing set on the isotoped annulus equivalent to
$T2_{2k-1}^+$.  Thus, $T2_{2k+1}^\pm$, $k \in {\bf \mathbb{Z}^+}$
can be reduced to $T2_1^\pm$. \endproof
\begin{lem} $T1_k$, $k \neq 0 \in {\bf \mathbb{Z}}$ can be reduced to $T1_{\pm 1}$.
\label{Decomp2T1kRED}
\end{lem}
\proof  Suppose there is a tight
contact structure on $M$ with $\Gamma_A = T1_{k}$, $k \leq -2$. We
show that there is a bypass along $A$ so that isotoping $A$ past
this bypasses produces a dividing set on the isotoped annulus
equivalent to $T1_{k+1}$.  In this way $T1_{k}$, $k \leq -2$ can be
reduced to $T1_{-1}$. The proof for positive $k$ is analogous.

\begin{figure}[ht!]
\vspace{0.2cm}
\begin{center}
\includegraphics[bb=0 0 116 139]{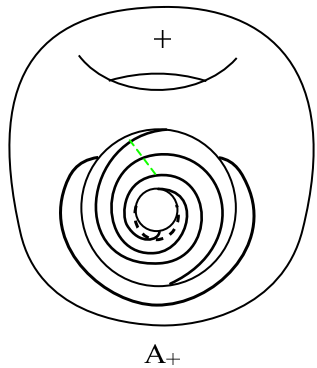}
\end{center}
\vspace{-0.5cm}
\caption{Convex decomposition for $M$ with $\Gamma_A$ = $T1_{-k}$}
\label{Decomp2T1minusRED}
\end{figure}

Consider the proposed bypass indicated in figure
\ref{Decomp2T1minusRED}. After rounding edges, we see that this
bypass is trivial. Although pictured for $T1_{-2}$, this is the
case for all $k \leq -2$. Isotoping $A_+$ across this bypass
yields a dividing set on the isotoped annulus equivalent to
$T1_{k+1}$. \endproof

Now we will turn our attention to constructing potentially tight
contact structures on $M$ by further decomposing (abstractly) $M_1
= M \backslash A$ with $\Gamma_A$ of type $T1_0$, $T1_1$,
$T1_{-1}$ and $T2_1^+$. For each of these four possible $\Gamma_A$
in turn, we will apply the gluing theorem
\cite{colgl,convexdecomp} or partition the constructed contact
structures into equivalence classes of tight and overtwisted
structures on $M_1$ by applying the gluing/classification theorem
for handlebodies. Further equivalences will be established on $M$
by locating state transitions along $A$. Gluing across the annulus
will be addressed in the following section.

\begin{figure}[ht!]
\begin{center}
\includegraphics[bb=0 0 266 255]{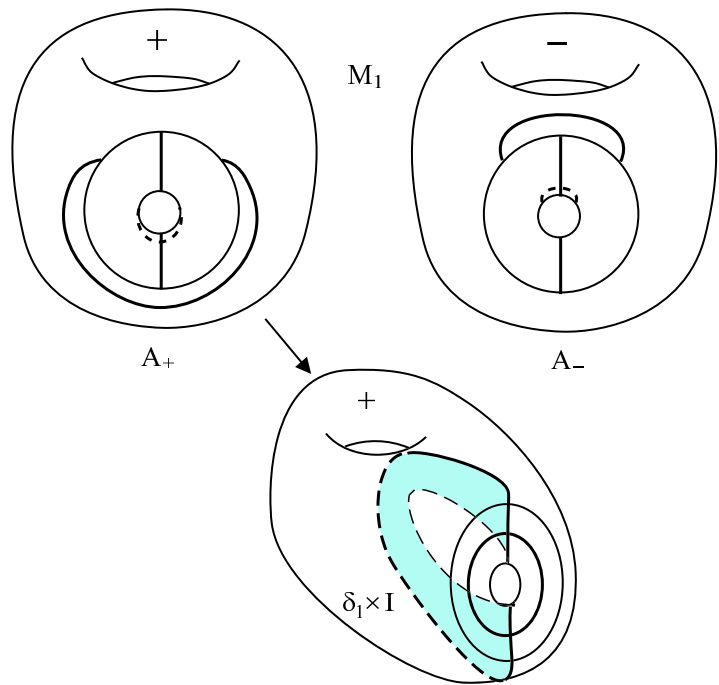}
\end{center}
\caption{Convex decomposition of $M$ with $\Gamma_A = T1_0$ (1)}
\label{Decomp2T10unique1}
\end{figure}

\begin{lem} There is a unique potentially tight contact structure on $M$ of
type $T1_0$, $T1_1$ and $T1_{-1}$.  These contact structures are all equivalent
on $M$ via state transitions. Moreover, they are universally tight
on $M_1 = M \backslash A$. \label{Decomp2T10POTTGT}
\label{T10}
\end{lem}

\proof \nopagebreak Consider $M_1 = M
\backslash A$ with $\Gamma_A = T1_0$ as pictured in figure
\ref{Decomp2T10unique1}. We show the decomposition for the component
of $M_1$ containing $A_+$ ($M_1^+$).  The argument for $M_1^-$ is
virtually identical.

\begin{figure}[ht!]
\begin{center}
\includegraphics[bb=0 0 180 91]{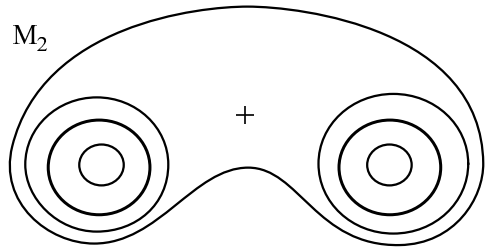}
\end{center}
\vspace{-0.3cm}
\caption{Convex decomposition of $M$ with $\Gamma_A = T1_0$ (2)}
\label{Decomp2T10unique2}
\end{figure}

\begin{figure}[b!]
\vspace{0.3cm}
\begin{center}
\includegraphics[bb=0 0 266 139]{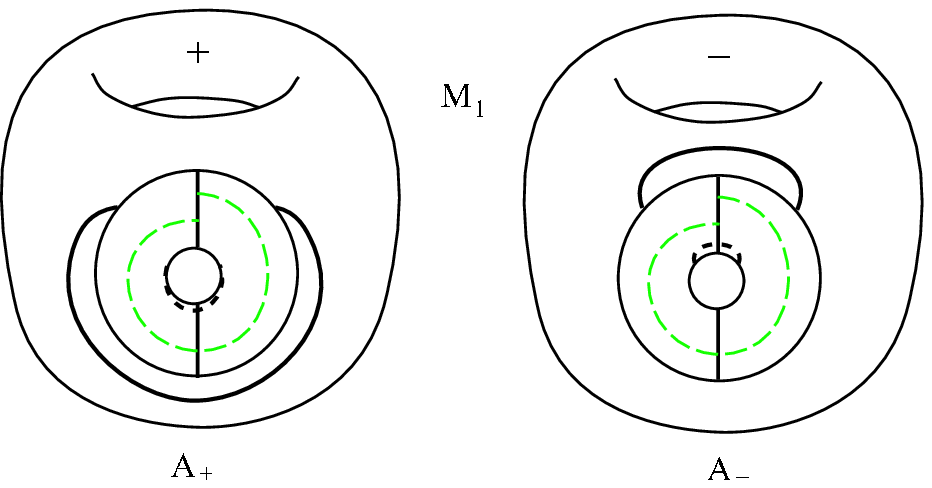}
\end{center}
\vspace{-0.5cm}
\caption{$M_1 = M \backslash A$} \label{M1isMminusAT102}
\end{figure}

We cut open $M_1^+$ along the convex cutting surface $\epsilon =
\delta_1 \times I$ with Legendrian boundary as indicated at the
bottom of figure \ref{Decomp2T10unique1}.  Two copies of the cutting
surface, $\epsilon_+$ and $\epsilon_-$ appear on the cut-open
manifold $M_2 = M_1 \backslash \epsilon$.  Since $tb(\partial
\epsilon_\pm) = -1$, there is only one tight possibility for
$\Gamma_\epsilon$.  Assuming this choice for $\Gamma_\epsilon$ and
rounding edges along $\partial \epsilon_\pm$, we get $(M_2 \cong S^1
\times D^2,\Gamma_{\partial_{M_2}})$ as in figure
\ref{Decomp2T10unique2}.

\begin{figure}[ht!]
\begin{center}
\includegraphics[bb=0 0 282 177]{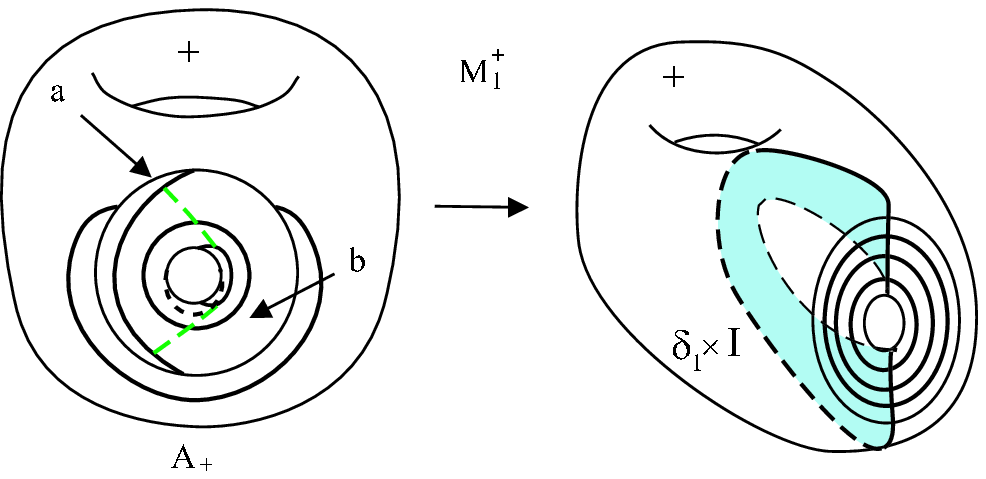}
\end{center}
\vspace{-0.3cm}
\caption{$M_1^+ = M_1 \backslash \epsilon$ with $\Gamma_{A_+} =
T2_1^+$ and possible bypasses} \label{M1T21possbyp}
\end{figure}

We further decompose $M_2$ by cutting along a convex, meridional
cutting surface $\eta$.  Since $tb(\partial \eta_\pm) = -1$, there
is only one tight possibility for $\Gamma_\eta$. Assuming this
choice for $\Gamma_\eta$ and rounding edges along $\partial
\eta_\pm$, we get $M_3 \cong B^3$ with  $\# \Gamma_{\partial_{M_3}}
= 1$. By Eliashberg's uniqueness theorem, there is a unique,
universally tight contact structure on $M_3 \cong B^3$ which extends
the one on the boundary.  Moreover, since the dividing sets on our
cutting surfaces $\epsilon$ and $\eta$ are boundary-parallel, we may
apply the gluing theorem \cite{colgl,convexdecomp} (theorem
\ref{GLUING}) to conclude the existence of a unique, universally
tight contact structure on $M_1$ with $\Gamma_A = T1_0$. If we
consider $M_1$ with $\Gamma_{A_+} = T1_{\pm 1}$ and round edges
along $\partial A_{\pm}$, we see that we obtain $M_1$ with a single
homotopically essential closed dividing curve as in the $T1_0$ case
(see figure \ref{Decomp2T10unique1}).  Proceeding with the convex
decomposition as in the $T1_0$ case shows that there is a unique,
universally tight contact structure on $M_1$ with $\Gamma_{A_+} =
T1_{1}$ and another with $\Gamma_{A_+} = T1_{-1}$.

\begin{figure}[ht!]
\begin{center}
\includegraphics[bb=0 0 213 219]{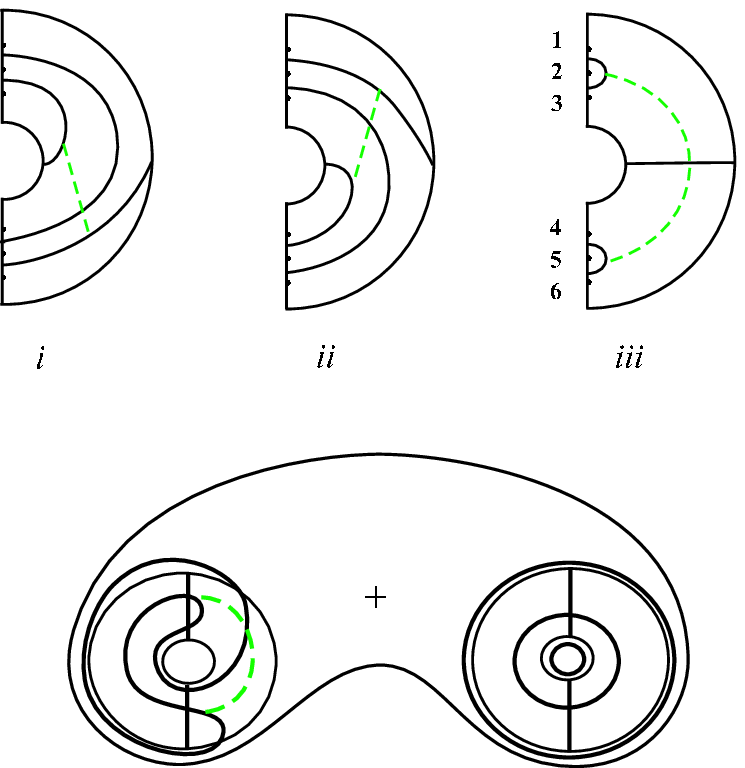}
\end{center}
\vspace{-0.3cm}
\caption{$M_2 = M_1 \backslash \epsilon$ with possible
$\Gamma_{\epsilon}$} \label{M1gammaetaposs}
\end{figure}

\begin{figure}[ht!]
\begin{center}
\includegraphics[bb=0 0 314 466]{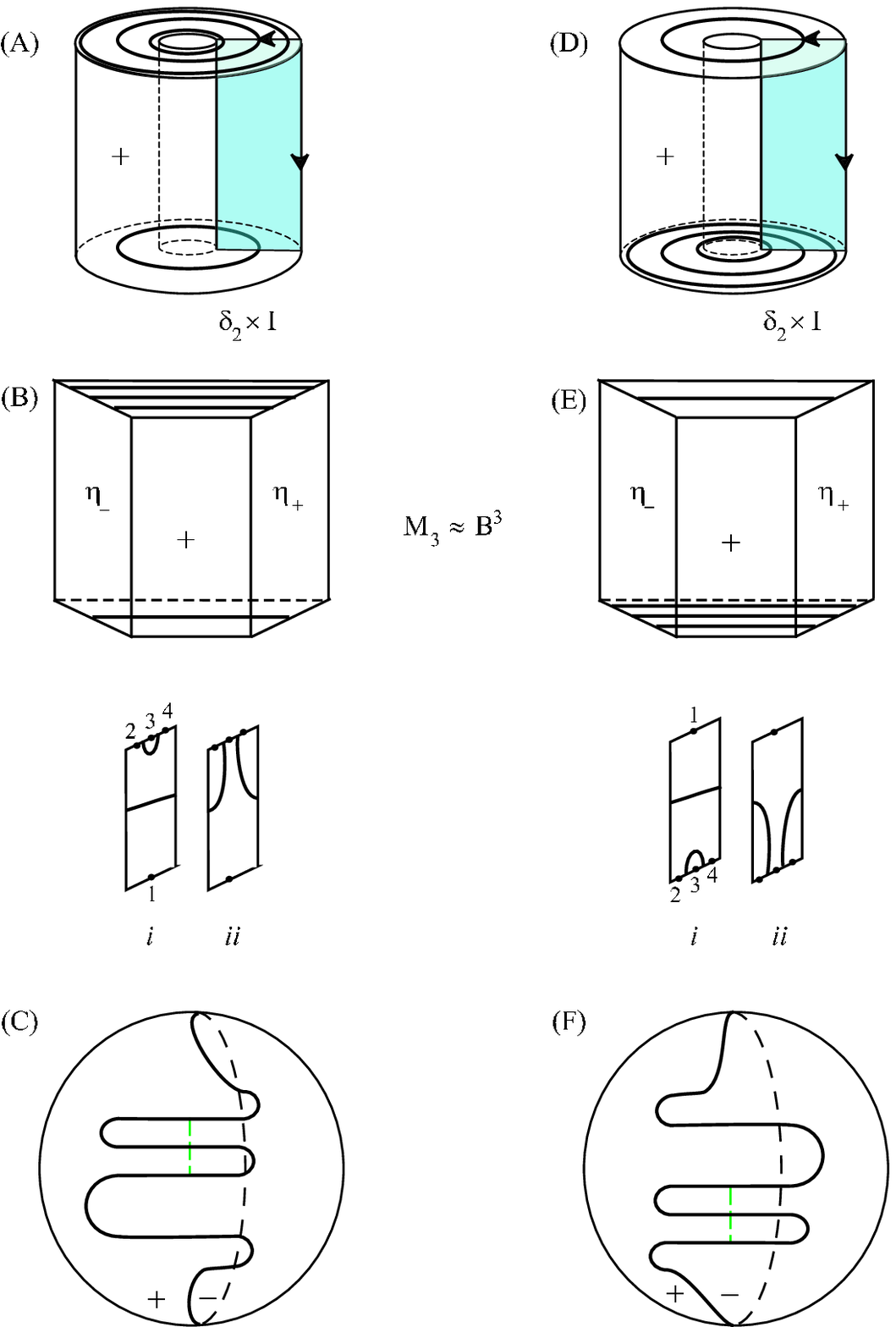}
\end{center}
\caption{Unique Non-Product $M_2 = M_1 \backslash \epsilon$ with
$\epsilon =  \emph{i}$ and \emph{ii}} \label{UniqueT21NP}
\end{figure}

Now we show the equivalence of the tight contact structures on $M_1$
with $\Gamma_{A_+} = T1_{0}$, $T1_{1}$ and $T1_{-1}$ by showing
there exists a state transition along $A$ transforming the unique,
universally tight contact structure on $M_1$ with $\Gamma_{A_+} =
T1_{0}$ into the unique universally tight contact structure on $M_1$
with $\Gamma_{A_+} = T1_{1}$ and another transforming it into the
unique universally tight contact structure on $M_1$ with
$\Gamma_{A_+} = T1_{-1}$. Let us consider $M_1$ with $\Gamma_{A_+} =
T1_0$ as in figure \ref{M1isMminusAT102}. The proposed bypasses
along $A_+$ and $A_-$ are trivial, and, hence, exist. They transform
$T1_0$ into $T1_1$ along $A_+$ and into $T1_{-1}$ along $A_-$ while
simultaneously transforming the dividing sets on $\epsilon$ and
$\eta$ for the decomposition in the $T1_0$ case into the dividing
sets on $\epsilon$ and $\eta$ for the decomposition in the $T1_1$
and $T1_{-1}$ cases, respectively. \endproof

By lemma 5.2 of \cite{FHM}, there is a unique universally tight
contact structure on $M \backslash A$ ($\amalg_{i = 1}^2 (S_i
\times I) $, $S_i \cong T^2 \backslash \nu(pt)$) with $\Gamma_A =
T1_0$, and it is given by perturbing the foliation of $M
\backslash A$ by leaves $S \times \{t\}$, $t \in [0,1]$. Recall
that we have fixed a characteristic foliation adapted to the
boundary-parallel dividing set $\Gamma_{\partial (M \backslash
A)}$.  The structure we describe here is the induced contact
structure in a product neighborhood of our disjoint union of
punctured tori such that $\mathcal{L}_{\frac{\partial}{\partial
t}}(\xi) = 0$. That is, flowing in the $I$- direction preserves
$\xi$.  This gives a dividing set on $A$ equivalent to $T1_0$. We
will call this ($I$- invariant) contact structure the
\emph{product structure}. Hence, lemma \ref{T10} states that the
unique, potentially tight structures on $M$ with $\Gamma_A =
T1_{\pm 1}$ are isotopic to the product structure.

\begin{lem} There is a unique, non-product potentially tight contact structure on $M$ of type $T2_1^+$.
This contact structure is tight on $M_1 = M \backslash A$.
\label{Decomp2T21POTTGT}
\end{lem}
\proof We will make the argument for the $A_+$
component of $M_1$ ($M_1^+$) since the argument for the other
component is completely analogous.

Suppose $\Gamma_A = T2_1^+$.  Note that if bypass $a$ indicated in
figure \ref{M1T21possbyp} exists, there is a decomposition of $M$
along an annulus $A^\prime$, isotopic to $A$, with
$\Gamma_{A^\prime} = T1_0$.  If bypass $b$ indicated in figure
\ref{M1T21possbyp} exists, there is a decomposition of $M$ along an
annulus $A^{\prime\prime}$, isotopic to $A$, with
$\Gamma_{A^{\prime\prime}} = T1_{-1}$. We first decompose $M_1$ in
order to isolate a potentially tight contact structure that does not
obviously contain one of these bypasses.

Let $\epsilon = \delta_1 \times I$ be the convex cutting surface
with Legendrian boundary indicated in figure \ref{M1T21possbyp}.
Cutting $M_1$ open along this cutting surface yields $M_2 = M_1
\backslash \epsilon$ with two copies $\epsilon_+$ and $\epsilon_-$
of the cutting surface.  All but the three choices of
$\Gamma_{\epsilon_+}$ given in figure \ref{M1gammaetaposs}
immediately lead to a homotopically trivial dividing curve, and,
hence, by \emph{Giroux's criterion}, to the existence of an
overtwisted disk.  The boundary-parallel dividing curves of choice
\emph{iii} may be realized along $A_+$ as bypasses $a$ and $b$
transforming $T2_1^+$ into $T1_0$ and $T1_{-1}$.  We want to show
that there is a unique potentially tight potentially non-product
structure on $M$ with $\Gamma_{\epsilon_+} = i$ and another with
$\Gamma_{\epsilon_+} = ii$, and that the state transition indicated
in choice \emph{i} of figure \ref{M1gammaetaposs} exists and takes
one to the other.  We will also show that there is no state
transition with bypass attaching arc indicated in choice \emph{ii}
of figure \ref{M1gammaetaposs} taking the unique potentially tight
potentially non-product structure with $\Gamma_{\epsilon_+} = ii$
into any structure with $\Gamma_{\epsilon_+} = iii$.  We then show
there is, in fact, a single potentially tight contact structure on
$M$ with $\Gamma_A = T2_1^+$ and $\Gamma_{\epsilon_+} = iii$
containing bypasses $a$ and $b$ as in figure \ref{M1T21possbyp}.
Finally we show that if the bypass indicated in figure
\ref{M1gammaetaposs} \emph{iii} exists inside $M_2$, then isotoping
$\epsilon$ past this bypass transitions this structure into an
obviously overtwisted one with $\Gamma_\epsilon = i$. The
state-transition argument for possible bypass digs along the
corresponding dividing sets for $\Gamma_{\epsilon_-}$ is similar.
This exhausts all potential states and state transitions on $M_1$
with $\Gamma_A = T2_1^+$ in this equivalence class under
gluing/classification (theorem \ref{gluingclassification}
\cite{gluing}).

\begin{figure}[ht!]
\begin{center}
\scalebox{.90}{\includegraphics[bb=0 0 406 447]{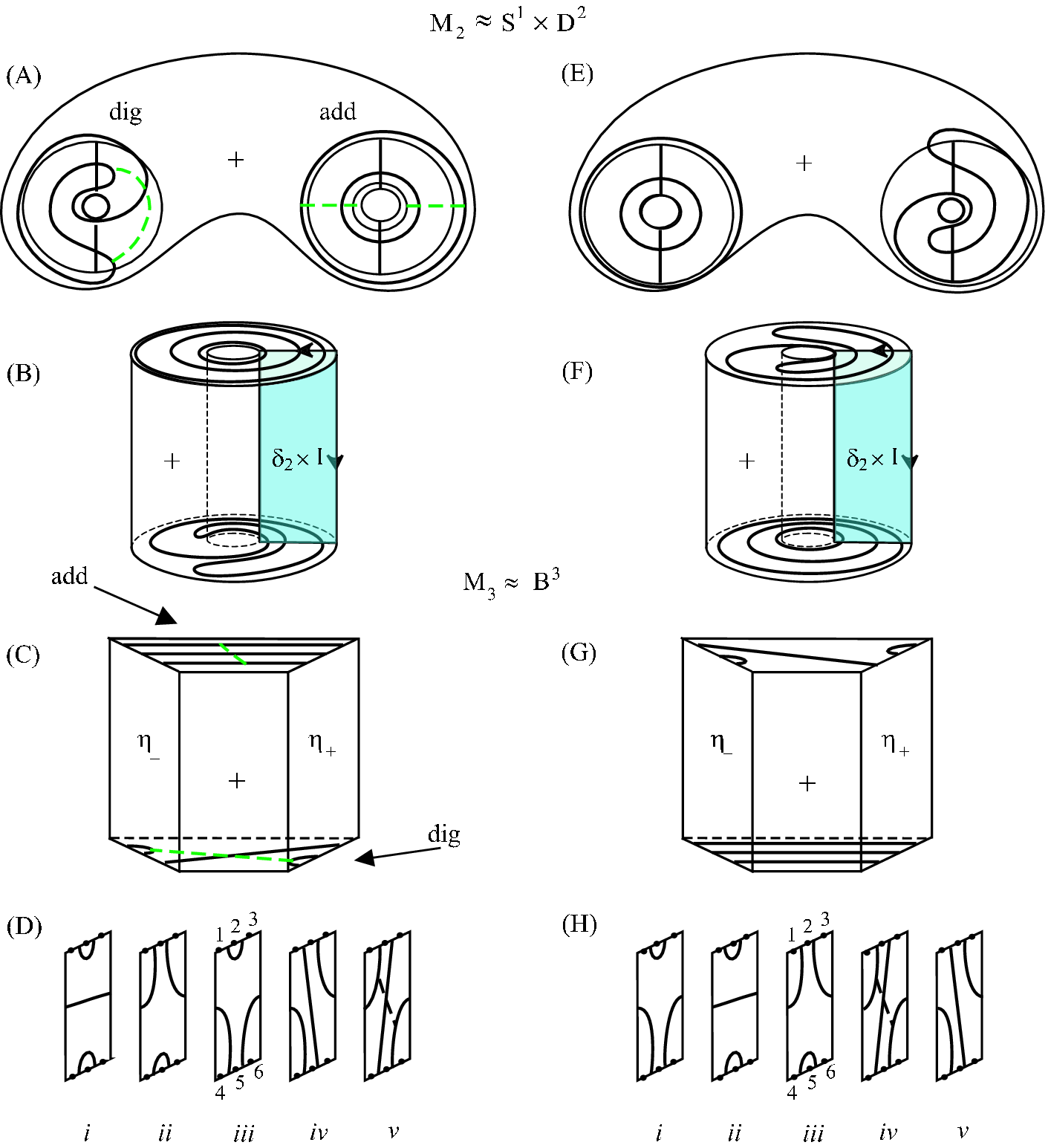}}
\end{center}
\caption{Equivalence of $M_2 = M_1 \backslash \epsilon$ with
$\epsilon =  \emph{i}$ and \emph{ii}} \label{UniqueT21NPST}
\end{figure}

Suppose we have $M_2$ with $\Gamma_{\epsilon_+} = i$. After rounding
edges, we get $M_2 \cong S^1 \times D^2$ with four longitudinal
dividing curves (see figure \ref{M1gammaetaposs}). Cutting $M_2$
open along a convex, meridional cutting surface $\eta = \delta_2
\times I$ as in figure \ref{UniqueT21NP} (A) yields $M_3 \cong B^3$
with two copies of the cutting surface, $\eta_+$ and $\eta_-$. Since
$tb(\partial \eta_+)$ = -2, there are two possibilities for
$\Gamma_{\eta_+}$ (see figure \ref{UniqueT21NP} (B)). One choice
contains a boundary-parallel dividing arc straddling position 3.
This indicates the existence of a bypass half-disk $B$ that can be
realized along $A_+$. Isotoping $A_+$ across $B$ produces a dividing
set on the isotoped annulus equivalent to $T1_0$. Applying the
remaining choice and rounding edges leads to $\# \Gamma_{\partial
B^3} = 1$ as in figure \ref{UniqueT21NP} (C). By Eliashberg's
uniqueness theorem, there is a unique extension of this contact
structure to the interior of $B^3$.  Since the dividing set on
$\eta$ is boundary-parallel, we may apply the gluing theorem
(theorem \ref{GLUING} \cite{colgl,convexdecomp}), we can conclude
that this contact structure is tight on $M_2$ with $\epsilon = i$.
If we similarly decompose $M_2$ with $\Gamma_{\epsilon_+} = ii$ as
in the figures \ref{UniqueT21NP} (D) through (F), we see that there
is a unique, potentially tight contact structure (that is not
obviously isotopic to a structure with a cut of type $\Gamma_A =
T1_0$) on $M$ that is tight on $M_2$ with $\Gamma_{\epsilon_+} =
ii$.

Now, consider the proposed state transition from choice \emph{i} to
choice \emph{ii} of $\Gamma_{\epsilon_+}$ on $M_2$.  The bypass
indicated on the left of figure \ref{UniqueT21NPST} (A)  is a
\emph{folding} bypass  on $M_2 \cong S^1 \times D^2$, and such
bypasses always exist \cite{convexdecomp}.  We need to show that if
we dig the bypass from $\epsilon_+$ and glue it back along
$\epsilon_-$, we transform the potentially tight potentially
non-product $T2_1^+$ with $\Gamma_{\epsilon_+} = i$ into the unique
potentially non-product $T2_1^+$ with $\Gamma_{\epsilon_+} = ii$. To
see this, we need to use a slightly different decomposition. Our new
decomposition for $M_2$ with $\Gamma_{\epsilon_+} = i$ appears in
figures \ref{UniqueT21NPST} (A) through (D).  Our new decomposition
for $M_2$ with $\Gamma_{\epsilon_+} = ii$ appears in figures
\ref{UniqueT21NPST} (E) through (H). In order to avoid cutting
through the proposed bypass attaching arc, our first cutting surface
(see figures \ref{UniqueT21NPST} (B) and (F)) will be a convex
cutting surface with Legendrian boundary that is not
\emph{efficient}.
\eject

In this decomposition, we first would like to isolate the
potentially tight potentially non-product contact structure from
the previous decomposition. By cutting $M_2$ inefficiently along
$\eta = \delta_2 \times I$, we have five choices for
$\Gamma_{\eta_+}$ as in figures \ref{UniqueT21NPST} (D) and (H).
Choices \emph{i} and \emph{ii} for both decompositions lead to a
homotopically trivial dividing curve and hence, by \emph{Giroux's
criterion}, an overtwisted disk.  In each decomposition, choice
\emph{iii} contains a boundary-parallel dividing curve (straddling
positions 2 and 5, respectively).  We conclude the existence of a
bypass half-disk that can be realized on $A_+$ in each
decomposition, transforming $\Gamma_{A_+}$ into $T1_0$ or
$T1_{-1}$. The structures represented by choices \emph{iv} and
\emph{v} are equivalent by the state transitions indicated in
figures \ref{UniqueT21NPST} (D) and (H). That is, (1) Both choices
\emph{iv} and \emph{v} in figure (D) yield, after edge rounding, a
single dividing curve on $\partial B^3$, (2) the abstract bypass
dig indicated in figure (D) \emph{v} transforms \emph{v} into
\emph{iv}, and (3) performing the bypass dig of figure (D)
\emph{v} does not change $\#\Gamma_{\partial B^3}$ (i.e.\ the
bypass disk exists inside the cut-open manifold). The same is true
for figures (H) \emph{iv} and \emph{v}. They represent the unique
potentially non-product $T2_1^+$ with $\Gamma_{\epsilon_+} = i$
and the unique potentially
non-product $T2_1^+$ with $\Gamma_{\epsilon_+} =
ii$, respectively (see figure \ref{UniqueT21NP}).

Now that we have identified the potentially tight potentially
non-product structure of figure \ref{UniqueT21NP} under this new
decomposition, we wish to show the existence of the state transition
along $\epsilon_+$ with $\Gamma_{\epsilon_+} = i$. Note that digging
the (folding) transitioning bypass from a portion of $\epsilon_+$ on
$M_3 \cong B^3$ and gluing it back along a portion of $\epsilon_-$
as in figure \ref{UniqueT21NPST} (C) transforms the unique
potentially tight potentially non-product structure on $M$ with
$\Gamma_{\epsilon_+} = i$ and $\eta_+ = iv$ into the unique
potentially tight potentially non-product structure on $M$ with
$\Gamma_{\epsilon_+} = ii$ and $\eta_+ = v$ .  This establishes
equivalence.

\begin{figure}[ht!]
\begin{center}
\includegraphics[bb=0 0 185 94]{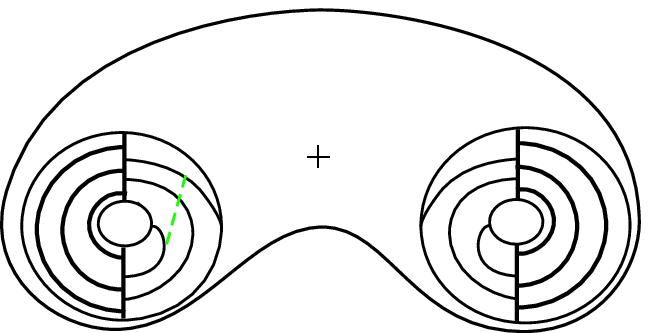}
\end{center}
\vspace{-0.3cm}
\caption{Non-Equivalence of $M_2 = M_1 \backslash \epsilon$ with
$\epsilon =  \emph{ii}$ and \emph{iii}} \label{Decomp2epsiloniii}
\end{figure}

We know from the previous decomposition that the potentially tight
potentially non-product contact structure on $M$ with
$\Gamma_{\epsilon_+} = i$ or $ii$ is tight on $M_2$. We now show
that the potential state transition taking $\Gamma_{\epsilon_+} =
ii$ into any structure with $\Gamma_{\epsilon_+} = iii$ does not
exist (see figures \ref{M1gammaetaposs} \emph{ii}, \ref{UniqueT21NP}
(F) and \ref{Decomp2epsiloniii}). From our choice of $\eta_+ = ii$,
we know that there is a bypass on the solid torus $M_2$ straddling
position 4 as in figure \ref{UniqueT21NP} (E).  This is equivalent
to adding a bypass straddling position 3 along the outside of the
torus (see the Attach=Dig property, p.64-66 of \cite{convexdecomp}).
Since both this add and the proposed transitioning bypass cannot
exist inside a tight manifold, we can conclude that the state
transition from $\Gamma_{\epsilon_+} = ii$ to $\Gamma_{\epsilon_+} =
iii$ does not exist.  Note that by this reasoning, we may conclude
that bypasses $a$ and $b$ in figure \ref{M1T21possbyp} cannot exist
in this case.

\begin{figure}[ht!]
\begin{center}
\includegraphics[bb=0 0 123 143]{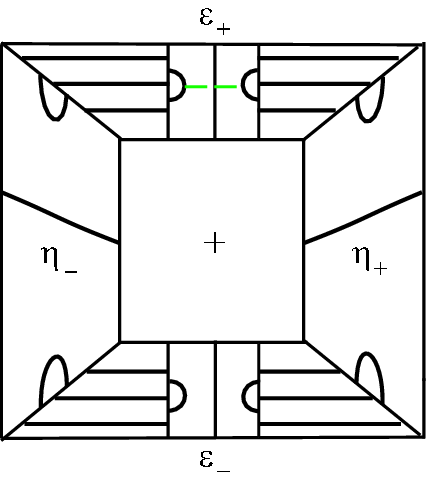}
\end{center}
\vspace{-0.3cm}\caption{Non-Equivalence of $M_2 = M_1 \backslash \epsilon$ with
$\epsilon =  \emph{iii}$ and \emph{i}} \label{Newiii}
\end{figure}

There is a single potentially tight structure on $M_1$ with
$\Gamma_A = T2_1^+$ and $\Gamma_{\epsilon_+} = iii$, and it is
pictured in figure \ref{Newiii}. All other potential $\Gamma_{\eta}$
are obviously overtwisted after edge-rounding. Clearly this choice
of $\Gamma_{\epsilon_+}$ indicates the existence of a bypass
half-disk that can be realized along $A$. Isotoping $A$ past this
bypass transforms $\Gamma_A = T2_1^+$ into $\Gamma_A = T1_0$.

Suppose that the bypass dig from $\epsilon_+$, as indicated in
figure \ref{Newiii}, exists (see also figure \ref{M1gammaetaposs}
\emph{iii}). From the figure, we see that digging the bypass from
$\epsilon_+$ and gluing back along $\epsilon_-$ transforms this
structure into an overtwisted one with $\Gamma_{\epsilon_+} = i$,
not into the potentially tight structure of figures
\ref{UniqueT21NP} (A) through (C) and \ref{UniqueT21NPST} (A)
through (D). The state-transition argument for possible bypass digs
along the corresponding dividing sets for $\Gamma_{\epsilon_-}$ is
similar.

Having checked all possible states and state transitions on $M_1$
with $\Gamma_A = T2_1^+$, we conclude, by Honda's
gluing/classification theorem \cite{gluing} (theorem \ref{GLUING}),
that the potentially tight potentially non-product structure on
$M_1$ with $\Gamma_{A_+} = T2_1^+$ is tight on $M_1$. This contact
structure is \emph{potentially allowable} on $M$.

\begin{figure}[ht!]
\begin{center}
\includegraphics[bb=0 0 236 136]{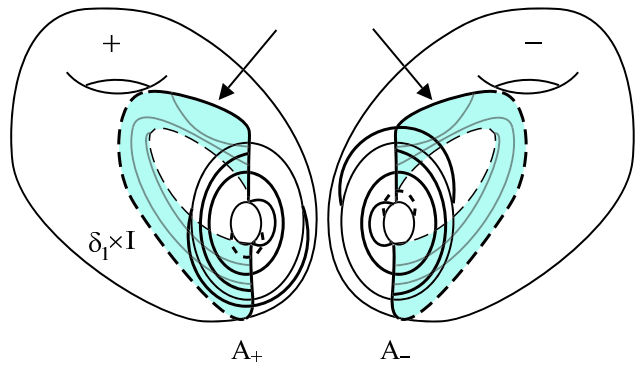}
\end{center}
\vspace{-0.3cm}
\caption{M contains a non-trivial bypass along $\Sigma_0$}
\label{NPbyp}
\end{figure}

\begin{figure}[ht!]
\begin{center}
\includegraphics[bb=0 0 203 127]{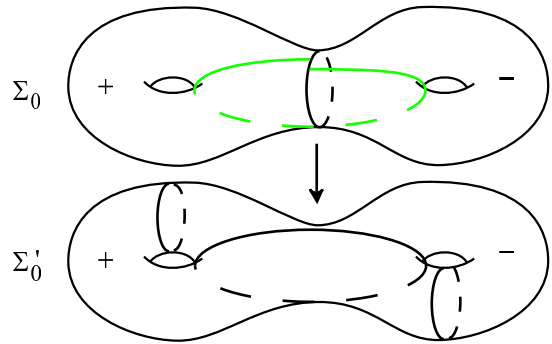}
\end{center}
\caption{Isotoping $\Sigma_0$ past a non-trivial bypass}
\label{NPbyp2}
\end{figure}

By the semi-local Thurston-Bennequin inequality (see
\cite{GirCircBund}, and proposition \ref{localtb} of this paper) we
know that there can be no non-trivial bypasses inside a contact
product neighborhood of a surface. Note that $M$ contains a
$\partial$-parallel arc as in figure \ref{NPbyp}, indicating the
existance of a bypass half-disk abutting $\partial M$. Isotoping
$\Sigma_0 = \Sigma_2 \times \{0\}$ past this bypass half-disk
results in a dividing set on the new, isotoped surface
$\Sigma_0^\prime$ as in figure \ref{NPbyp2}. Since this non-trivial
bypass exists along $\partial M$, we may conclude that this
structure is, in fact, not the product structure. \endproof

We now turn to the task of showing that this contact structure comes
from a unique, non-product tight contact structure on $M$.

\section{Gluing}
It is now necessary to establish the tightness of the two
potentially allowable (non-product and product) contact structures
of the previous section. The strategy here will be to use the ideas
involved in the proof of the gluing theorem for contact manifolds
with convex boundary and boundary-parallel dividing curves on all
gluing surfaces. This proof was originally given by Colin
\cite{colgl} and subsequently formulated in terms of convex
decompositions by Honda et al.~\cite{convexdecomp}.  We say the
dividing set $\Gamma_S$ on a convex surface is
\emph{boundary-parallel} if $\Gamma_S$ is a collection of arcs
connecting $\partial S$ to $\partial S$ and this collection of arcs
cuts off disjoint half-disks along the boundary of $S$. A contact
structure $\xi$ on $M$ is \emph{universally tight} if the pull-back
$(\tilde{M}, \tilde{\xi})$ of the contact structure to any cover of
$M$ is tight.  The general statement 
was given earlier (the ``gluing theorem'' theorem \ref{GLUING}).

There are two main obstacles to applying the gluing theorem
directly.  First of all, the dividing set $T2_1^+$ on $A$ is not
boundary-parallel.  So, we cannot use the theorem directly to
$M_1$ glued along $A$. In the gluing theorem, the
boundary-parallel requirement is necessary in order to guarantee
that any bypass along the gluing surface \emph{at most} introduces
a pair of parallel dividing curves. However, we know from our
decomposition precisely which bypasses exist along $A$. In the
non-product case, they are trivial bypasses and folding bypasses
along the central homotopically non-trivial closed dividing curve
of $T2_1^+$, which introduce a pair of dividing curves parallel to
the original curve.

The second obstacle is that it is necessary in our case to use a
\emph{state transition} argument in establishing tightness of the
non-product contact structure on $M_1  = M \backslash A$.  This
argument relies on Honda's gluing/classification theorem which
guarantees tightness but not \emph{universal} tightness.  Thus, we
cannot pull our contact structure back to an arbitrary cover and
expect that the structure remains tight.  Instead, we will
construct explicit covers $\tilde{M_i}$ of $M \backslash A$ and
compute pull-back structures directly in order to establish
tightness of $(\tilde{M_i},\tilde{\xi})$. Therefore, it is
possible to establish the conditions necessary to apply the ideas
of the gluing theorem and conclude tightness of $(M,\xi)$.

The idea of the proof here, following the proof of the gluing
theorem, will be as follows.  First, we will construct finite covers
of $M_1$ in which all of the aforementioned folding bypasses are
trivial and prove that these covers with the pull-back contact
structures are tight.  Then, we will assume the existence of an
overtwisted disk $D$ inside M and look at controlled pull-backs of
the bypasses necessary to push $A$ off of $D$ to the specified tight
covers. In this manner, we will construct a cover
$(\tilde{M},\tilde{\xi})$, a pull-back of $\tilde{A}$ of $A$ and a
lift $\tilde{D}$ of the overtwisted disk $D$.  In this cover, all
the bypasses needed to isotope $\tilde{A}$ off of $\tilde{D}$ will
be trivial.  Using tightness of $(\tilde{M},\tilde{\xi})$, we can
derive a contradiction to the existence of $D$, thereby establishing
tightness of $M$.
\subsection{Constructing tight covers}
\begin{figure}[ht!]
\begin{center}
\includegraphics[bb=0 0 258 223]{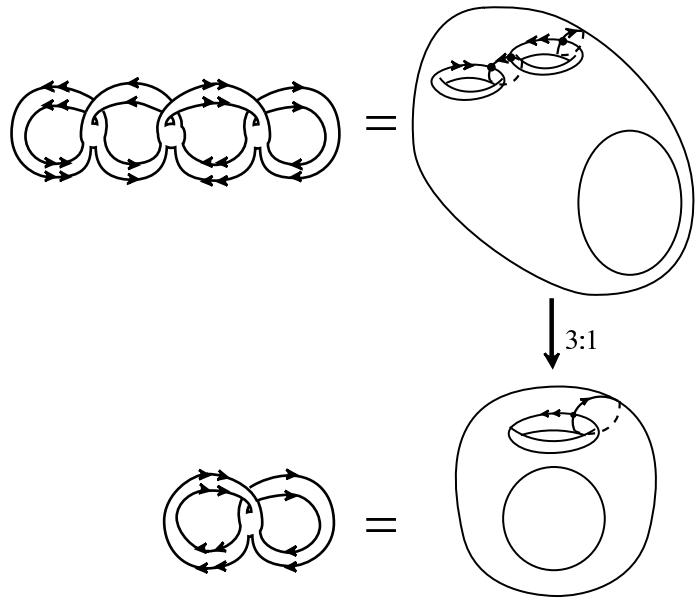}
\end{center}
\caption{A covering space for $T^2 - \nu(pt)$} \label{3to1cover}
\end{figure}

Let us begin by constructing a 3:1 cover $\tilde{S}$ of the
punctured torus as in figure \ref{3to1cover}. Since this cover has
a single boundary component, it must be a once punctured surface
$\Sigma_g$ for some $g \in \mathbb{Z}^+$. An Euler characteristic
calculation tells us that this cover $\tilde{S}$ must be $\Sigma_2
- \nu(pt)$:
\begin{eqnarray*}
  \chi(\tilde{S}) &=& 1-2g \\
  \empty &=& 3(1-2(1))\\
  \empty &=& -3
\end{eqnarray*}

\begin{figure}[ht!]
\begin{center}
\includegraphics[bb=0 0 134 57]{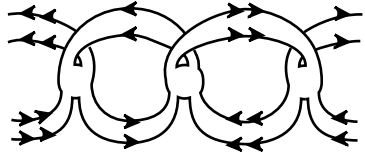}
\end{center}
\caption{Fundamental domain for the 3:1 cover of $(\Sigma_2 -
\nu(pt)) \times I$} \label{funddom}
\end{figure}

Thickening each surface by crossing with an interval induces a
cover of the punctured torus cross $I$ by $(\Sigma_2 - \nu(pt))
\times I$. Thus, we have constructed a 3:1 cover of each component
$M_1 = M \backslash A$.  Now, we use this construction and the
fundamental domain in figure \ref{funddom} to construct a 3:1
cover of $(\Sigma_2 - \nu(pt)) \times I$.  In this way, we
construct a $3^n$\negmedspace:1 cover of $M_1$ for each $n \in
\mathbb{Z}^+$.

\begin{figure}[ht!]
\begin{center}
\includegraphics[bb=0 0 316 224]{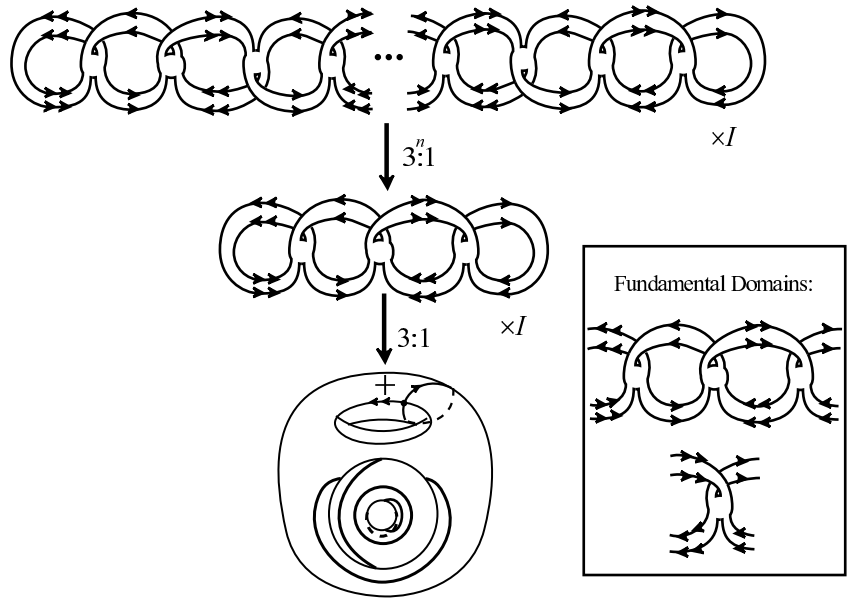}
\end{center}
\caption{A covering space for $M_1$} \label{coveringspace}
\end{figure}

We now have a $m$:1-fold covering space $(\tilde{M}=\Sigma_{m+1}
\times I,p)$ of $M=\Sigma_2 \times I$ such that $m = 3^n$ for each
$n \in \mathbb{Z}^+$.  The restriction of such a cover to $M_1$ is
the disjoint union of two copies of $(\Sigma_{\frac{m+1}{2}} -
\nu(pt)) \times I$ (see figure \ref{coveringspace}).  Note that
the lift $\tilde{A}$ of the annulus $A$ is another annulus
(\emph{enlarged} from the original by a factor of $m$).  Let
$M_1^+$ denote the $A_+$ component of $M \backslash A$. We will
establish tightness of $(\tilde{M} \backslash \tilde{A}
,\tilde{\xi})$ where $\xi$ is the unique, non-product tight
contact structure on $M \backslash A$ by focusing on the covering
space $(\tilde{M} \backslash \tilde{A},p|_{p^{-1}(M_1^+)})$. The
argument for the other component is completely analogous.  Recall
that the product structure is universally tight on $M_1$.

\begin{figure}[ht!]
\begin{center}
\includegraphics[bb=0 0 321 282]{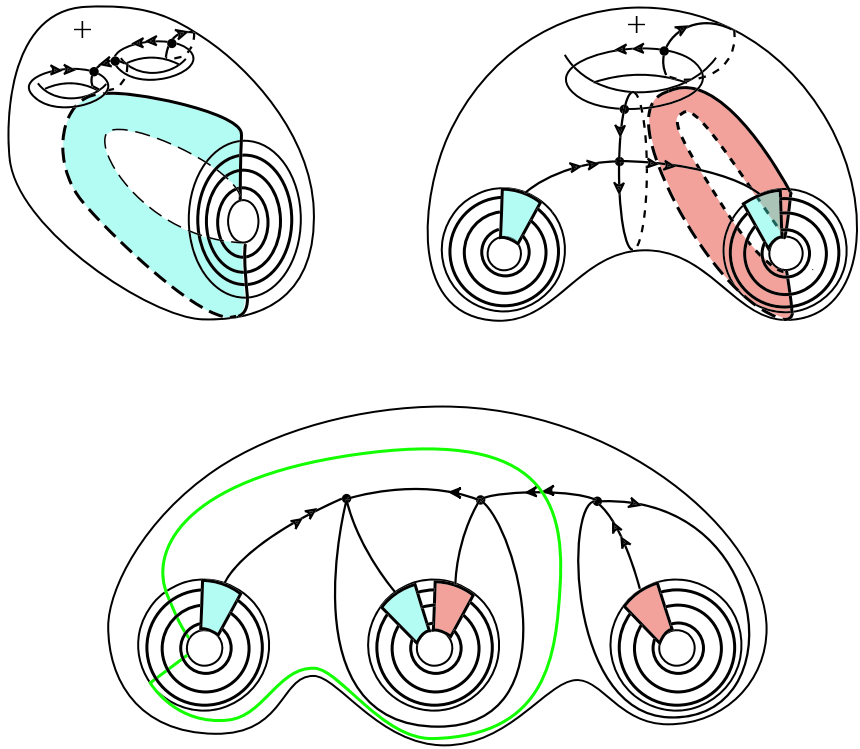}
\end{center}
\caption{The pullback of the contact structure}
\label{cutsupstairs1}
\end{figure}

\begin{lem} $(\tilde{M} \backslash \tilde{A} ,\tilde{\xi})$ is a tight contact manifold
where $\xi$ is the unique, non-product tight contact structure on $M \backslash A$.
\label{PullbackIsTight}
\end{lem}
\proof  It suffices to consider
$(\tilde{M} \backslash \tilde{A},p|_{p^{-1}(M_1^+)})$.  Our aim is
to prove tightness of $(\tilde{M} \backslash \tilde{A}
,\tilde{\xi})$ by using Honda's gluing/classification theorem on
the convex decomposition of $(\tilde{M} \backslash \tilde{A}
,\tilde{\xi})$ which is the pullback of the one on $M \backslash
A$ (see lemma \ref{Decomp2T21POTTGT} and figures \ref{UniqueT21NP}
and \ref{UniqueT21NPST}). We will begin by lifting the cut on $M
\backslash A$ which is transverse to the double-arrow direction as
indicated for the 3:1 cover in figure \ref{cutsupstairs1}. Recall
that the 3:1 cover of $M_1^+$ is $\Sigma_2 - \nu(pt)$.

\begin{figure}[ht!]
\begin{center}
\includegraphics[bb=0 0 325 233]{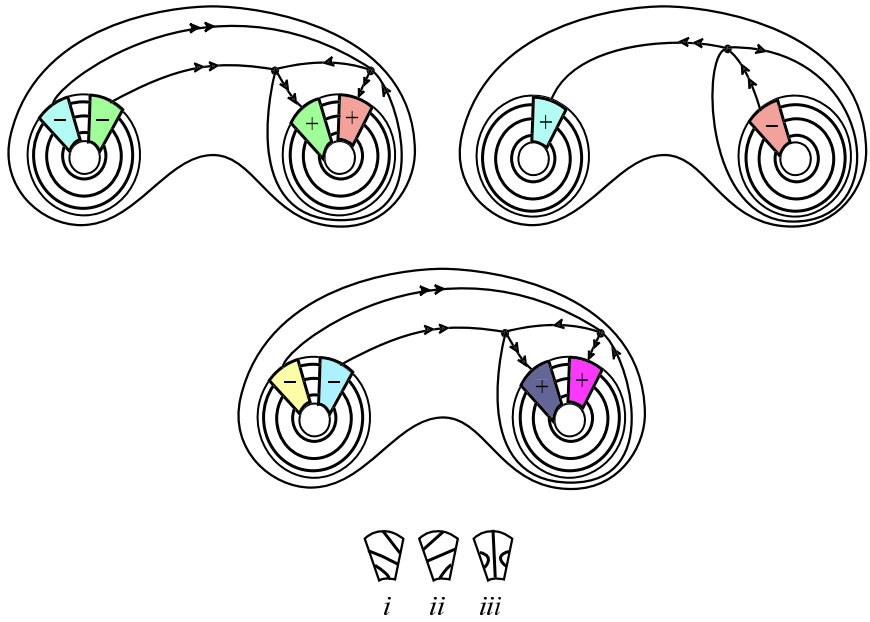}
\end{center}
\vspace{-0.3cm}
\caption{The pullback of the contact structure $\# 2$}
\label{cutsupstairs2}
\end{figure}

 Cutting open along the $m$ pull-backs of the first cutting disk
downstairs yields $\frac{m-1}{2} + 1$ solid tori.  Our covering
projection $p$ restricted to $\frac{m-1}{2}$ of these tori is a
2:1 cover, while $p$ restricted to the remaining torus is a 1:1
cover. These $\frac{m-1}{2} + 1$ solid tori fall into three
categories according to the way the cutting surfaces appear on
them. One of the 2:1 covers always contains the positive and
negative copies of one cutting surface, the positive copy of a
second cutting surface, and the negative copy of a third. The 1:1
cover always contains the positive copy of a cutting surface and
the negative copy of another.  The remaining $\frac{m-1}{2} - 1$
2:1 covers contain copies of four different cutting surfaces, two
positive and two negative. Examples of these three categories are
given in figure \ref{cutsupstairs2} (note that no tori of the
latter form appear in the decomposition of the 3:1 cover).

\begin{figure}[ht!]
\vspace{0.3cm}
\begin{center}
\includegraphics[bb=0 0 380 205]{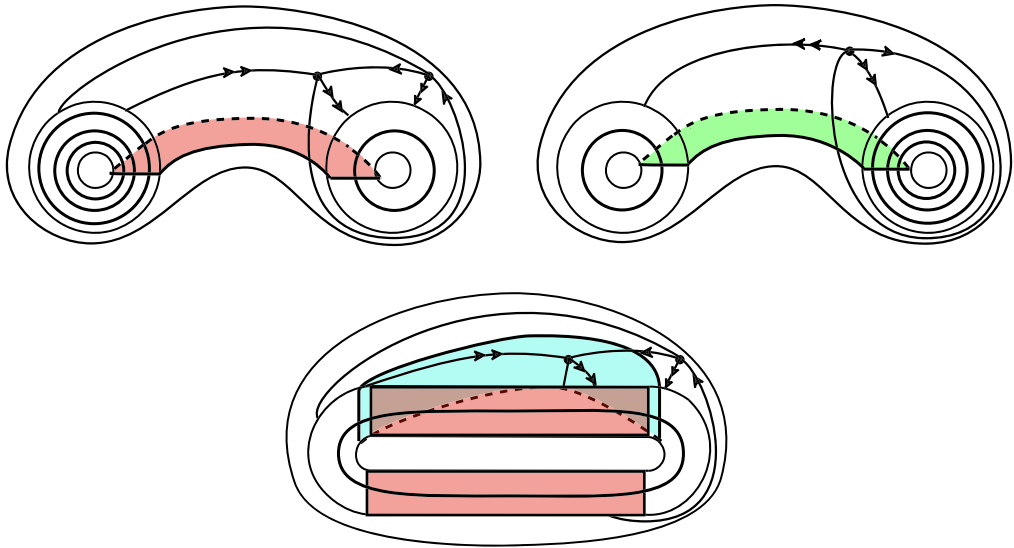}
\end{center}
\vspace{-0.3cm}
\caption{The pullback of the contact structure $\# 3$}
\label{cutsupstairs3}
\end{figure}

Now we begin pulling back the contact structure.  At this stage,
this means applying dividing curve configuration \emph{i} to the
cutting surfaces (see lemma \ref{Decomp2T21POTTGT} and figure
\ref{M1gammaetaposs}).  The resulting boundary configurations
(parallel sets of longitudinal dividing curves on the tori) are
given for the 3:1 cover at the top of figure \ref{cutsupstairs3}
along with three new convex cutting disks which are the pull-backs
of the second cut downstairs (transverse to the single arrow
direction). Pulling back the contact structure at this level means
applying the dividing curve configuration of figure
\ref{UniqueT21NP} (B) \emph{ii} to these disks. The result is $m$
3-balls with a single dividing curve each.

\begin{figure}[ht!]
\begin{center}
\includegraphics[bb=0 0 311 138]{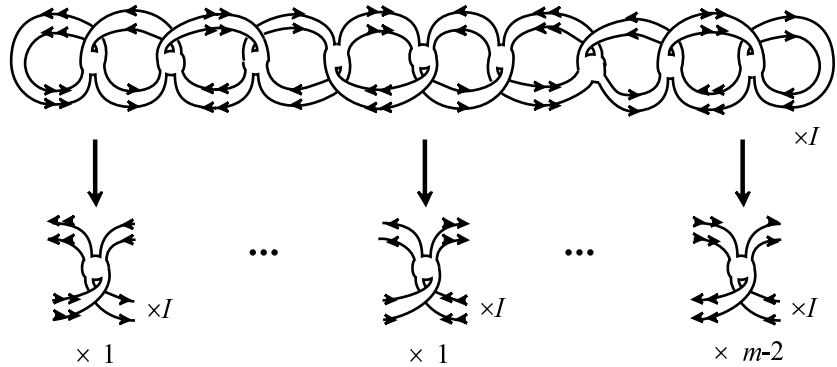}
\end{center}
\vspace{-0.3cm}
\caption{Decomposition of $\tilde{M}$ into three types of balls}
\label{checksonballs}
\end{figure}

This state is not obviously overtwisted (i.e.\ it is
\emph{potentially allowable}).  Recall that the state transition $i$
to $ii$ along $\epsilon_+$ exists downstairs while the state
transition $ii$ to $iii$ does not. The state transition $i$ to $ii$
downstairs corresponds to doing all such transitions along all lifts
of $\epsilon_+$ upstairs.  It is necessary to check that doing any
combination of these bypasses upstairs transitions us to a
potentially allowable state. These state transitions can be checked
explicitly and exist as trivial or folding bypasses on each torus.
The results are dividing curve configurations on the tori consisting
of either two, four or six parallel, longitudinal dividing curves.
The non-existence of the state transitions $ii$ to $iii$ and $i$ to
$iii$ downstairs does not imply the non-existence of any such state
transition upstairs.  However, it can be checked that this
possibility never exists.  The finiteness of the check is a result
of the fact that, for any cover, the pull-back decomposition results
in the union of $m = 3^n$ balls (each containing two copies of the
pullback of each of the two cutting disks downstairs).  These balls
fall into three categories. Two of the $m$ balls contain two
different types of self-gluing, while the remaining $m-2$ balls
admit no self-gluing (see figure \ref{checksonballs}).   Thus, we
may conclude, by the gluing/classfication theorem, that $(\tilde{M}
\backslash \tilde{A} ,\tilde{\xi})$ is tight. \endproof

\subsection{An interpretation of the gluing theorem and tightness}
In order to establish tightness of the product structure we simply
note that the pull-back $\tilde{\xi}$ of the product structure is a
product structure, and this structure is tight by Giroux's
criterion.  For the potentially allowable non-product structure, we
note that the proof of the standard gluing theorem (theorem
\ref{GLUING} \cite{convexdecomp}) depends on two facts that we don't
have:
\begin{enumerate}
    \item The dividing sets on cutting surfaces are
    boundary-parallel.
    \item The contact structure $\xi$ is universally tight on $M \backslash S$ where $S$ is our cutting
    surface.
\end{enumerate}

However, we have classified our contact structures on $M \backslash
A$ and know which bypasses are possible.  By reviewing the proof of
the gluing theorem, we see that the requirement that $\Gamma_S$ be
boundary-parallel is there so that the types of bypasses possible
along $S$ are strictly limited to those which are trivial or
``long''. In our case the types of bypasses possible along $A$ with
$\Gamma_A = T2_1^+$ inside $(M,\xi)$ are similarly limited (where
$\xi$ is the unique non-product structure). We have the following
lemma.

\begin{figure}[ht!]
\begin{center}
\includegraphics[bb=0 0 146 66]{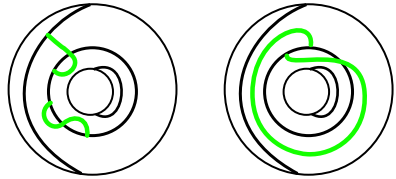}
\end{center}
\vspace{-0.3cm}
\caption{Some bypass possibilities} \label{BypPoss}
\end{figure}

\begin{lem} Let $A$ be a convex annulus whose boundary is as indicated in figure
\ref{MisMtimesI2} with $\Gamma_A = T2_1^+$ inside the unique
non-product $(M,\xi)$. Then, any convex annulus $A^\prime$ that is
obtained from $A$ by a sequence of bypass moves will have a dividing
set $\Gamma_{A^\prime}$ differing from $\Gamma_A$ by possibly adding
an even number of parallel, closed dividing curves encircling the
inner boundary component $\delta \times 1$ of $A$.
\label{PossBypasses}
\end{lem}
\proof Recall that, for the
non-product structure on $M \backslash A$ with $\Gamma_A = T2_1^+$,
the bypasses indicated along $A_+$ at the left of figure
\ref{M1T21possbyp} and their counterparts along $A_-$ do not exist.
These bypasses also do not exist in any of the covers constructed in
the previous section for similar reasons (see figure
\ref{cutsupstairs3} and lemma \ref{Decomp2T21POTTGT}). By examining
all remaining possibilities, we see that the only bypasses that
exist along $A$ in this case are trivial bypasses and folding
bypasses.  Further bypasses may be trivial (i.e.\ they produce no
change in the dividing set), may add pairs of parallel dividing
curves, or may delete pairs of dividing curves.  Examples of these
possibilities are given in figure \ref{BypPoss}. So, the dividing
set on any annulus $A^\prime$ obtained from $A$ by a sequence of
bypass moves has dividing set $\Gamma_{A^\prime}$ that differs from
$A$ by adding, at most, pairs of homotopically non-trivial closed
dividing curves parallel to the original one. \endproof

To address our second problem, let us consider the strategy of the
of the proof of the gluing theorem. The proof is by contradiction.
The existence of an overtwisted disk $D \subset M$ is supposed, and
we consider the sequence of bypass moves required to isotope $A$ off
of $D$. The requirement that $\xi$ is universally tight on $M
\backslash A$ is so that, when analyzing a single bypass move along
$A$ (which is trivial or increases $\# \Gamma_A$), we can lift to a
large enough cover so as to make this bypass trivial. We continue to
lift to covers until we arrive at the cover $\tilde{M}$ of $M$ with
lifts $\tilde{A}$ of $A$ and a lift $\tilde{D}$ of a proposed
overtwisted disk $D$ in which all bypasses needed to isotope
$\tilde{A}$ off of $\tilde{D}$ (producing the isotoped annulus
$\tilde{A^\prime}$) are trivial. But $\Gamma_{\tilde{A^\prime}}$ =
$\Gamma_{\tilde{A}}$ where $\tilde{A^\prime}$ is obtained from
$\tilde{A}$ by a sequence of trivial bypass attachments. Since $\xi$
is universally tight on $M \backslash A$, this gives us a
contradiction.

The objective of the previous section was to construct specific
covers of $\tilde{M}$ so that $\tilde{M} \backslash \tilde{A}$ is
tight. These covers suffice to satisfy the requirements of the
gluing theorem.  Thus, we are able to establish tightness of the
potentially tight non-product $(M , \xi)$.

We repeat here two lemmas from the proof of the gluing theorem
\cite{convexdecomp}. The first one concerns isotoping $A$ off of
an overtwisted disk $D$.  By a result of Honda \cite{gluing}, we
may perturb the characteristic foliation on the cutting surface
and look at the local model near points on the boundary of the
disk. We can arrange for $D$ to be transverse to $A$ and for
$\partial D$ to be contained in $\Gamma_A$.  Moreover, after
possibly modifying $D$, we can assume that the hypotheses of the
Legendrian realization principle are satisfied so that $D \cap A$
consists of Legendrian arcs and curves. This is called the
``\emph{controlled intersection}'' of an overtwisted disk with a
convex cutting surface $A$.

Now, to push $A$ away from $D$ so that we eliminate a closed curve
of intersection $\delta$.  Let $D_\delta$ denote the subdisk of $D$
with $\partial D_\delta = \delta$. Since we assume $A$ has a tight
neighborhood, we must have that $t(\delta,Fr_A) < 0$.  We may
perturb  $D_\delta$ (rel boundary) so as to make it convex with
Legendrian boundary.

We have that  $\stackrel{\circ}D_\delta$ is contained in $M
\backslash A$, which is tight.  So we may assume that the dividing
set of $D_\delta$ consists only of embedded arcs with endpoints on
$\partial D$.  We may push $A$ to engulf the bypass which
corresponds to one of these boundary-parallel dividing curves on
$D_\delta$. We continue until all bypasses are consumed.

For an arc of intersection $\gamma$ in $D \cap A$, we proceed
similarly, but we choose a disk $D_\gamma$ with $\gamma \subset
\partial D_\gamma$.  This concludes the sketch of the following
lemma:

\begin{lem}  It is possible to isotope $A$ off of $D$ in a finite number of steps, each of which is a bypass along $A$.
\label{FiniteBYP}
\end{lem}

The second lemma concerns the effect of isotoping cutting surfaces
across trivial bypasses.  A trivial bypass along a convex surface
$A$ can be realized inside an $I$-invariant (``contact product'')
neighborhood of the surface.

\begin{lem}  If $A$ is a convex surface with Legendrian boundary inside a contact manifold $(M , \xi)$ and $A^\prime$ is a convex surface obtained from $A$ by a trivial bypass, then $A$ and $A^\prime$ are contact isotopic and, hence, $(M \backslash A, \xi)$ is tight if and only if $(M \backslash A^\prime, \xi)$ is tight.
\label{TrivialBYP}
\end{lem}

We are now ready to prove an adaption of the gluing theorem. Much of
this proof is identical to the general proof.

\begin{thm} The potentially tight non-product structure is tight on $M$.
\label{Decomp2GLUING}
\end{thm}

\proof Assume that $M$ is not
tight. Then, there exists an overtwisted disk $D \subset M$. We
can perform a contact isotopy so that $D$ and $A$ intersect
transversely along Legendrian curves and arcs and so that
$\partial D \cap A \subset \Gamma_A$ \cite{gluing}. We note that
closed curves in $D \cap A$ are homotopically trivial in $A$.  We
want to eliminate the innermost closed curves on $D$ by pushing
$A$ across $D$.  Consider a two-sphere $S$ formed by a disk on $D$
and one on $A$ whose common boundary is an innermost curve of
intersection $\delta \subset D$. Then, $S$ bounds a ball across
which we can isotope A.

By lemma \ref{FiniteBYP}, we can push $A$ across $D$ to decrease the
number of intersections of $A$ with $D$ in a finite number of bypass
steps that possibly change the dividing curve configuration on the
isotoped annulus.  If we consider a single bypass along $A$, we know
it is either trivial or increases $\Gamma_A$ (see the proof of lemma
\ref{PossBypasses}).  Recall that the covering spaces constructed
in the last section enlarge $A$ by a factor of $m$. So we can lift
to a large enough cover of this type so that any folding bypass
attachment becomes trivial. We can continue lifting through covers
of this type until we arrive at the cover $\tilde{M}$ of $M$ with
lifts $\tilde{A}$ of $A$ and $\tilde{D}$ of a proposed overtwisted
disk $D$ in which all bypasses needed to isotope $\tilde{A}$ off of
$\tilde{D}$ are trivial. But $\Gamma_{\tilde{A^\prime}}$ =
$\Gamma_{\tilde{A}}$ where $\tilde{A^\prime}$ has $\tilde{A^\prime}
\cap \tilde{D} = \emptyset$ and is obtained from $\tilde{A}$ by a
sequence of trivial bypass attachments. Since $\tilde{\xi}$ is tight
on $M \backslash \tilde{A}$, then, by lemma \ref{TrivialBYP},
$\tilde{\xi}$ is tight on $M \backslash \tilde{A^\prime}$. This
contradicts the existence of an overtwisted disk in $M \backslash
\tilde{A^\prime}$.  Thus, $(M,\xi)$ is tight. \endproof

\section{A special property of the non-product tight
contact structure on $M$}

By viewing the non-product tight contact structure on $M = \Sigma
\times I$ from the perspective of a different convex decomposition,
it is possible to see that this contact structure contains every
possible bypass abutting the boundary.

\begin{figure}[ht!]
\begin{center}
\includegraphics[bb=0 0 200 146]{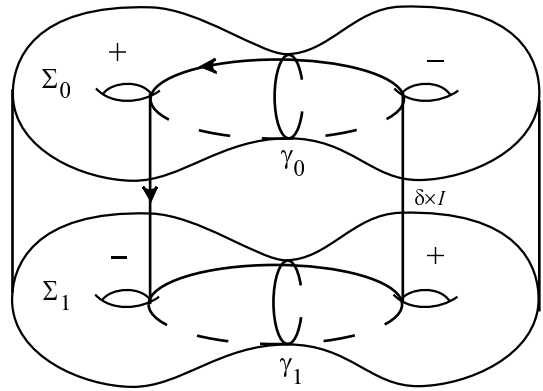}
\end{center}
\vspace{-0.3cm}
\caption{M = $\Sigma_2 \times I$} \label{MisMtimesI}
\end{figure}

\begin{figure}[b!]
\vspace{0.3cm}
\begin{center}
\includegraphics[bb=0 0 187 153]{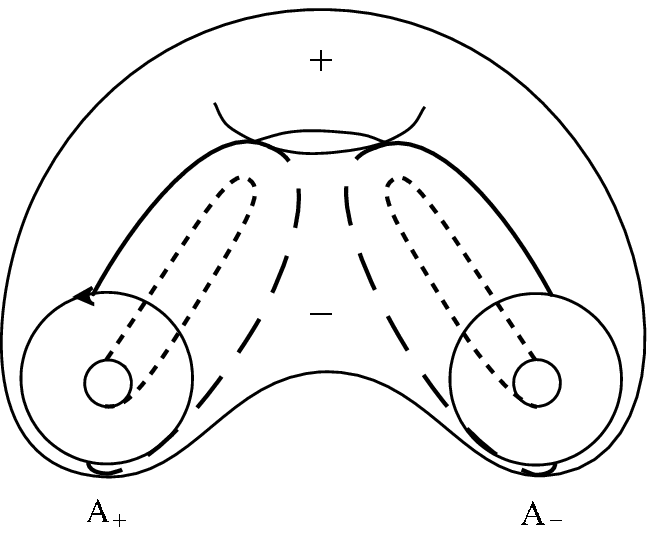}
\vspace{-0.3cm}
\caption{\label{M1isMminusA} $M_1=M \backslash A$}
\end{center}
\end{figure}

We perform a second convex decomposition assuming we are starting
with the non-product structure $\xi$ on $M$.  The first cut of
this decomposition will be along $\delta \times I$ where $\delta
\times \{i\}$ is a homotopically nontrivial non-separating curve
which is efficient with respect to $\Gamma_{\Sigma \times \{i\}}$,
$i \in [0,1]$ (see figure \ref{MisMtimesI}).  After this initial
annular cut, we have $M_1 = M \backslash A$ which is a twice
punctured torus cross $I$ as in figure \ref{M1isMminusA}. We have
the same four possibilities for $\Gamma_A$ as in the previous
decomposition (see figure \ref{GammaAtypes}), but in this case it
is easy to see that all $T2_{2m}^\pm$, $m \in \mathbb{Z}^+$ are
overtwisted.  We have the following lemmas. The proofs will be
omitted since they are similar to the proofs of the reduction
lemmas for the previous decomposition.
\begin{enumerate}
    \item Lemma:  $T2_{2m+1}^\pm, m \in \mathbb{Z}^+$ can be reduced to $T2_1^\pm$.
    \item Lemma:  $T1_k, k \in \mathbb{Z^+}$ can be reduced to $T1_1$.
    \item Lemma:  $T1_k, k \in \mathbb{Z^-}$ can be reduced to $T1_{-1}$.
\end{enumerate}

Since we have already shown that the non-product structure is tight,
we will focus our attention on showing that the non-product
structure, viewed from the perspective of the current decomposition,
has a cut of type $T2_1^+$ and one of type $T2_1^-$.  We will do
this by showing that there is a non-product contact structure in
each of the categories $T2_1^+$, $T2_1^-$, $T1_{1}$ and $T1_{-1}$
and we can find isotopies transforming each of these structures into
one-another.

\begin{thm} There exists exactly one universally tight contact structure on $M_1$, not equivalent to
$T1_{0}$, in each of the categories $T1_{\pm 1}$ and $T2_1^{\pm}$.
\label{theoremT2T1equiv}
\end{thm}

\begin{figure}[ht!]
\begin{center}
\scalebox{0.75}{\includegraphics[bb=0 0 442 608]{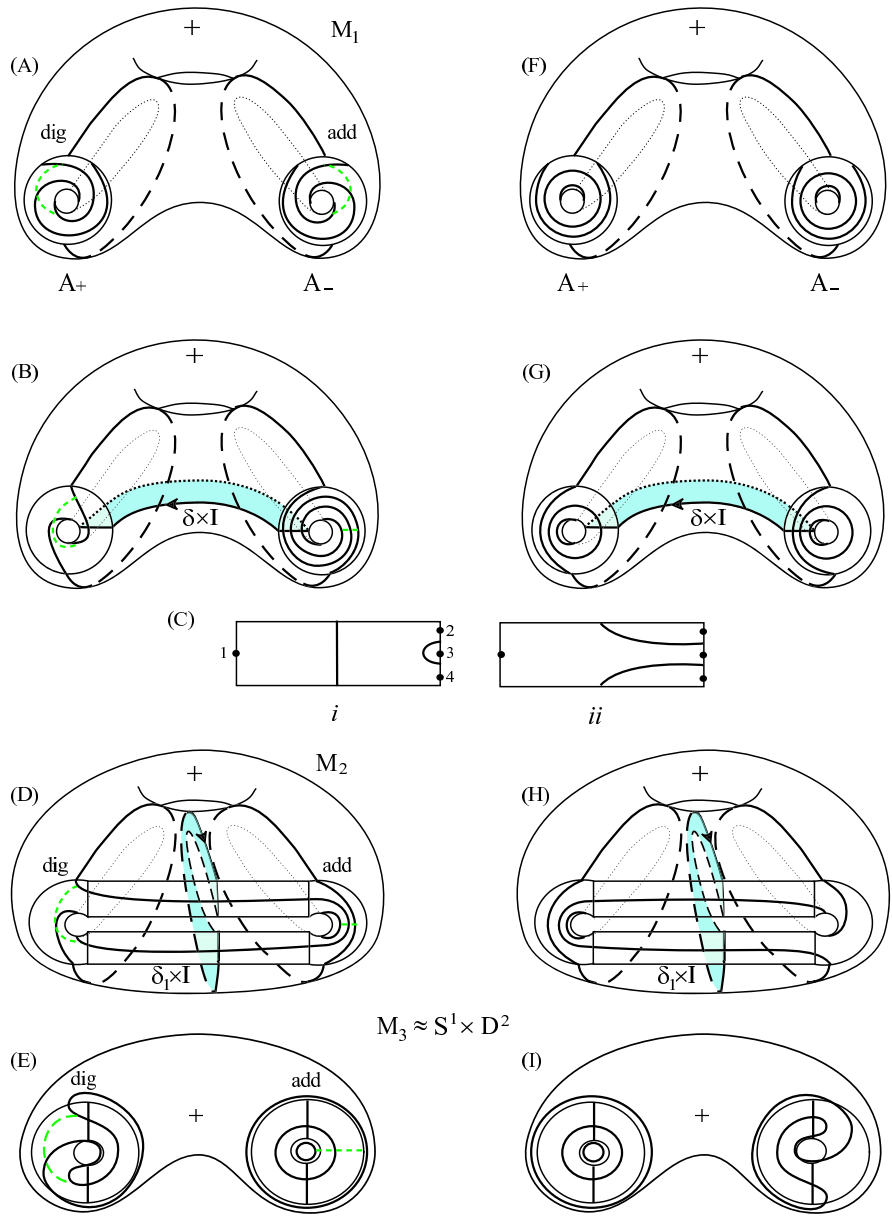}}
\end{center}
\caption{\label{T11T21plusunique1} Convex decomposition $\# 1$ for
$T1_1$ and $T2_1^+$ (1)}
\end{figure}

\begin{figure}[ht!]
\begin{center}
\includegraphics[bb=0 0 363 238]{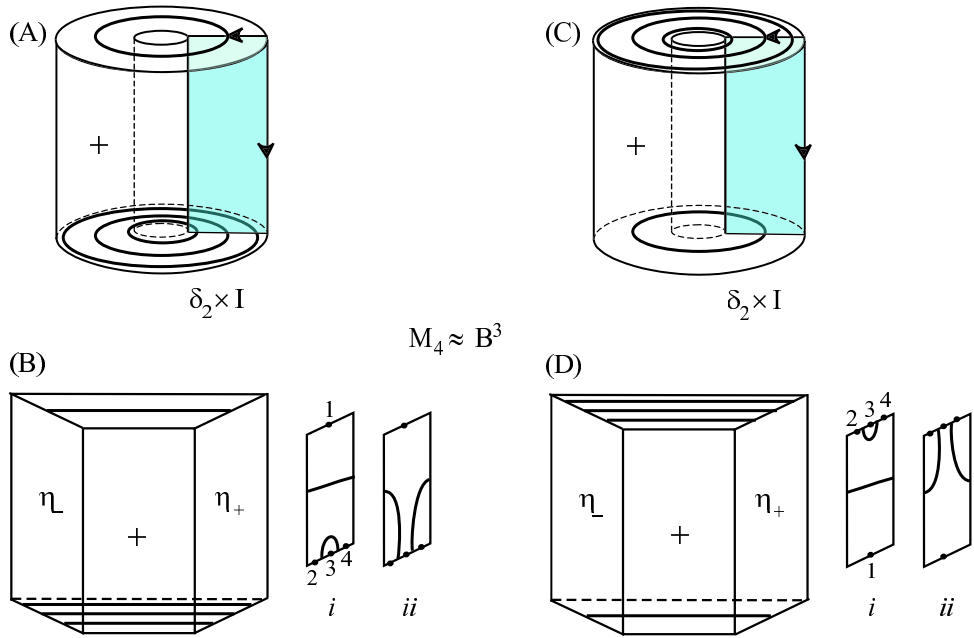}
\end{center}
\caption{\label{T11T21plusunique2} Convex decomposition $\# 1$ for
$T1_1$ and $T2_1^+$ (2)}
\end{figure}

\proof We show that there is
exactly one non-product tight contact structure on $M_1$ with
$\Gamma_A$ of type $T1_1$.    To do this, we provide a convex
decomposition of $M=\Sigma \times I$ starting with the convex
annulus $A$ together with $\Gamma_A=T1_1$.  The convex disks
defining this decomposition are given in figures
\ref{T11T21plusunique1} and \ref{T11T21plusunique2}.  The other
cases are argued similarly.

 \emph{Round} \emph{edges} along $\partial A_{\pm} \subset M_1$
and choose a new cutting surface $\delta \times I$  where $\delta$
is indicated in figure \ref{T11T21plusunique1} (B).  Assume $\delta
\times I$ is convex with \emph{efficient} Legendrian boundary. After
cutting along $\delta$, we have a new manifold $M_2$ which is a
punctured torus cross $I$ (see figure \ref{T11T21plusunique1} (D)).
Two copies of the cutting surface $\delta \times I$ appear in $M_2$
that we will call $\epsilon_+$ and $\epsilon_-$.  Since $tb(\partial
\epsilon_\pm) = -2$, there exist two possible dividing curve
configurations for $\Gamma_{\epsilon_+}$.  One of these
configurations induces a bypass half-disk straddling position 3 as
shown in figure \ref{T11T21plusunique1} (C) \emph{i}. Isotoping $A_-
\subset M_1$ across this bypass produces a new dividing set on the
isotoped annulus $A_-^\prime$ equivalent to $T1_0$.

Let $\Gamma_{\epsilon_+}$ be the remaining choice.  After
\emph{rounding} \emph{edges} along $\partial \epsilon_\pm$, we
obtain $(M_2,\Gamma_{\partial M_2})$ as in figure
\ref{T11T21plusunique1} (D).  After choosing a new convex cutting
surface $\delta_1 \times I$ with \emph{efficient} Legendrian
boundary as in figure \ref{T11T21plusunique1} (D), we cut open $M_2$
and obtain $(M_3 \cong S^1 \times D^2,\Gamma_{\partial M_2})$ with
two copies of the cutting surface $\delta_1 \times I$. Call these
copies $\tau_+$ and $\tau_-$.  Since $tb(\partial \tau_\pm) = -1$,
there is exactly one possibility for $\Gamma_{\tau_+}$.  Applying
this configuration to $\tau_+$ and $\tau_-$ and \emph{rounding}
\emph{edges} along $\partial \tau_\pm$ leads to the configuration
$(M_3,\Gamma_{\partial M_3})$ as shown in figure
\ref{T11T21plusunique1} (E).

We view $M_3$ as a solid torus and choose a convex, meridional
cutting surface $\delta_2 \times I$ with \emph{efficient} Legendrian
boundary as in figure \ref{T11T21plusunique2} (A). We cut open $M_2$
along this surface to obtain $(M_4 \cong B^3,\Gamma_{\partial M_2})$
with two copies of the cutting surface $\delta_1 \times I$.  Call
these copies $\eta_+$ and $\eta_-$. Since $tb(\partial \eta_\pm) =
-2$, there are two possibilities for $\Gamma_{\eta_+}$. One of these
configurations induces a bypass half-disk straddling position 3 as
shown in figure \ref{T11T21plusunique2} (B) \emph{i}. Isotoping $A_-
\subset M_1$ across this bypass produces a new dividing set on the
isotoped annulus $A_-^\prime$ equivalent to $T1_0$.

Let $\Gamma_{\eta_+}$ be the remaining choice.  After
\emph{rounding} \emph{edges} along $\partial \eta_\pm$, we obtain
$M_4 \cong B^3$ with $\#(\Gamma_{\partial M_4}) = 1$. Thus there is
at most one non-product tight contact structure on $M_1$ with
$\Gamma_A$ of type $T1_1$.

\begin{figure} [ht!]
\begin{center}
\scalebox{0.80}{\includegraphics[bb=0 0 433 486]{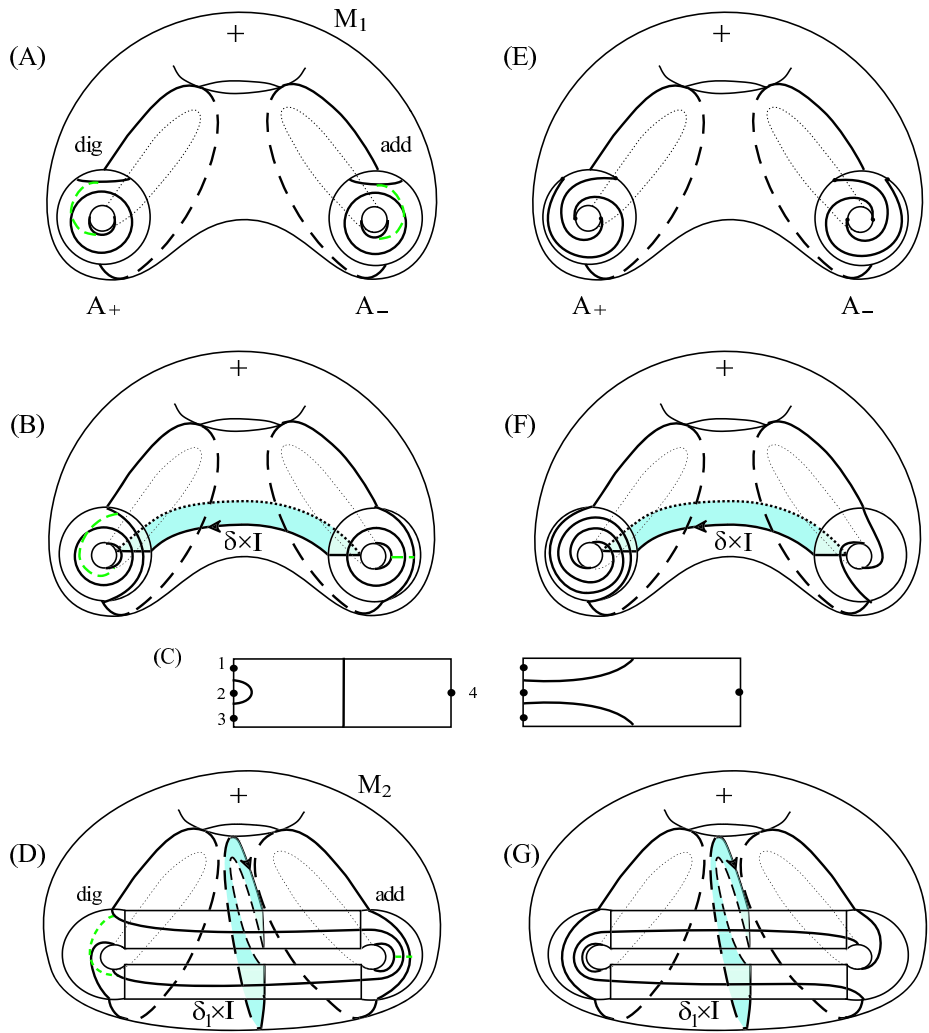}}
\end{center}
\caption{\label{T21minusT1-1unique1} Convex decomposition $\#$ 1 for
$T2_1^-$ and $T1_{-1}$}
\end{figure}

 Now, by Eliashberg's uniqueness theorem, there is a
unique \emph{universally} tight contact structure on $M_4 = B^3$
which extends the one on the boundary. Also, the dividing sets on
the convex disks $\epsilon$, $\tau$ and $\eta$ corresponding to our
unique non-product contact structure were all boundary parallel.
Therefore, by the gluing theorem \cite{colgl,FHM}, there is a
unique, universally tight contact structure on $M_1$ with $\Gamma_A$
of type $T1_1$. Similarly, we can establish the existence of a such
a structure structure on $M_1$ with $\Gamma_A$ of types $T2_1^+$
(given in figures \ref{T11T21plusunique1} (F) through (I) and (C)
and \ref{T11T21plusunique2}(C) and (D)), $T2_1^-$ (figure
\ref{T21minusT1-1unique1} (A) through (D))and $T1_{-1}$ (figure
\ref{T21minusT1-1unique1} (E) through (G) and (C)) .\endproof

\begin{lem} The unique, non-product universally tight contact structures on $M_1$ of type $T1_{\pm 1}$ and
$T2_1^{\pm}$ are all equivalent.
\end{lem}

\begin{figure} [ht!]
\begin{center}
\includegraphics[bb=0 0 335 502, width=0.85\hsize]{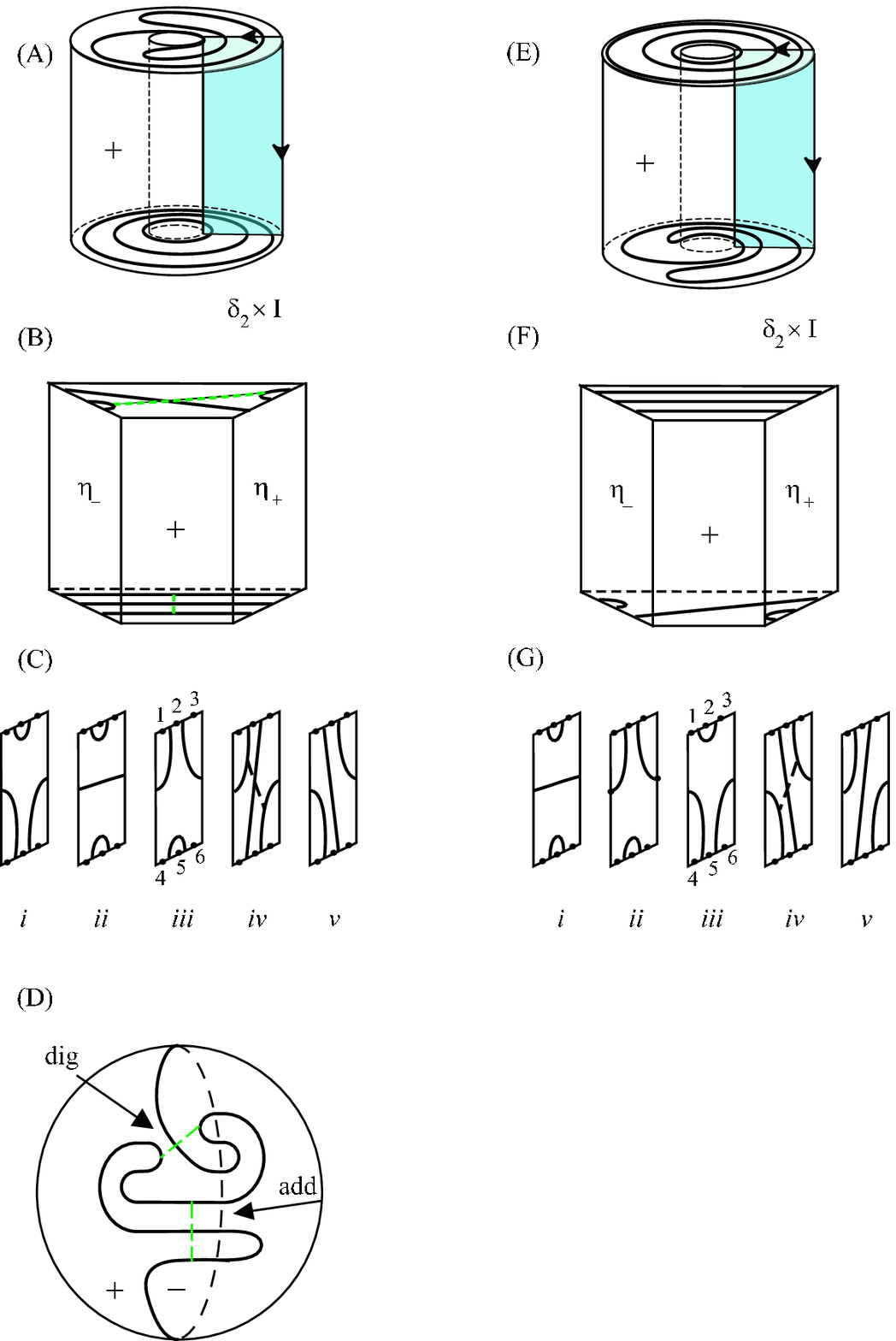}
\end{center}
\caption{\label{T11T21plusST} Equivalence of $T1_1$ and $T2_1^+$}
\end{figure}

\proof We prove the equivalences $T1_1 \cong
T2_1^+$.  The others are done similarly. We start by taking a convex
decomposition of $M=\Sigma \times I$ starting with the convex
annulus $A$ together with $\Gamma_A=T1_1$.  The convex disks
defining this decomposition are given in figures
\ref{T11T21plusunique1} and \ref{T11T21plusST}.  The purpose is to
find a bypass half disk $B$ as indicated in figure
\ref{T11T21plusunique1} (A) so that digging $B$ from one side of $A$
and adding it to the other transforms $\Gamma_A = T1_1$ into
$\Gamma_A = T2_1^+$. The other cases are argued similarly.

 First, consider the partial convex decomposition $M_3 =
M_1 \backslash (\epsilon \cup \tau)$ of the unique non-product
$T1_1$ as in the previous lemma and figures \ref{T11T21plusunique1}
(A) (B) (C) (D) and (E).  We now proceed with a slightly different
decomposition. View $M_3$ as a solid torus and choose a convex,
meridional cutting surface $\delta_2 \times I$ with Legendrian
boundary. In order to prove existence of the proposed bypass
half-disk $B$, our choice of cutting surface shown in figure
\ref{T11T21plusST} (A) is not \emph{efficient}.

We observe that the proposed bypass is a \emph{folding}
\emph{bypass}, and such bypasses always exist \cite{convexdecomp}.
However, it is easy to show existence explicitly in our case by
proceeding with the decomposition.

After cutting $M_3$ open along $\delta_2 \times I$, we obtain the
ball $(M_4,\Gamma_{\partial M_3})$.  Two copies of the cutting
surface $\delta_2 \times I$ appear in $M_4$.   Call them $\eta_+$
and $\eta_-$. There exist five choices for $\Gamma_{\eta_+}$ as
pictured in figure \ref{T11T21plusST} (C). Applying choice \emph{i}
or \emph{ii} of $\Gamma_{\eta_{+}}$ and \emph{rounding} \emph{edges}
along $\partial \eta_\pm$ yields a dividing set on $\partial M_4
\approx S^2$ consisting of three dividing curves. By \emph{Giroux's}
\emph{criterion}, we conclude the existence of an overtwisted disk.
Applying choice \emph{iii} of $\Gamma_{\eta_{+}}$ induces a bypass
half-disk straddling position 5. Isotoping  $A_- \subset M_1$ across
this bypass produces a new dividing set on the isotoped annulus
$A_-^\prime$ equivalent to $T1_0$.

Both choices \emph{iv} and \emph{v} correspond to a dividing curve
configuration on $S^2 \approx \partial B^3$ such that
$\#\Gamma_{\partial B^3}= 1$ ($\Gamma_{\partial B^3}$ resulting
choice \emph{iv} is illustrated in figure \ref{T11T21plusST} (D)).
We would like to find a  \emph{state} \emph{transition} taking
dividing curve configuration \emph{iv} to dividing curve
configuration \emph{v}. To do this, we must establish the
existence of a bypass half-disk as indicated by it's attaching arc
in choice \emph{iv} of figure \ref{T11T21plusST} (C). Applying
choice \emph{iv} to $\partial M_4 \approx S^2$, we see that this
bypass is the trivial one on the ball. Such a bypass is guaranteed
to exist by the \emph{right to life} principle \cite{FHM}. Thus,
we have shown the existence the proposed bypass $B$ taking $T1_1$
to $T2_1^+$.

We proceed by performing the same convex decomposition on $M_1$
with $\Gamma_A$ of type $T2_1^+$.  We then argue that digging $B$
yields a \emph{state} \emph{transition} transforming the unique
potentially non-product $T1_1$ into the unique potentially
non-product $T2_1^+$. The convex disks defining this decomposition
are given in figures \ref{T11T21plusunique1} (F) (G) (C) (H) and
(I) and \ref{T11T21plusST} (E) (F) and (G).

We proceed exactly as in the previous case, noting that, in the
third stage of the decomposition of $M_1$, the boundary of our
cutting surface $\eta$ shown in figure \ref{T11T21plusST} (E) is
not efficient. Since $tb(\partial \eta_\pm) = -3$, we will have
five possibilities for $\Gamma_{\eta_+}$.  Cases \emph{i} and
\emph{ii} in figure \ref{T11T21plusST} (G) lead to
$\#(\Gamma_{\partial B^3}) = 3$, and case \emph{iii} induces a
bypass that can be realized on $A_- \subset M_1$ transforming
$T2_1^+$ into $T1_0$.  Finally, there exists a \emph{state}
\emph{transition}, as indicated in figure \ref{T11T21plusST} (G)
\emph{iv} taking case \emph{iv} to case \emph{v}.

We will now show that, at each stage of the convex decomposition,
digging the bypass along a section of $A_+$ and adding it back along
a section of $A_-$ transforms the unique non-product $T1_1$ into the
unique non-product $T2_1^+$.  First, consider $M_1=M \backslash A$
with $\Gamma_A=T1_1$ as pictured in figure \ref{T11T21plusunique1}
(A).  If we dig the bypass along $A_+$ and add it back along $A_-$
as indicated, the result is exactly $M_1=M \backslash A$ with
$\Gamma_A=T2_1^+$ as pictured in figure \ref{T11T21plusunique1} (F).
Second, consider $M_2=M_1 \backslash (\delta \times I)$ as pictured
in figure \ref{T11T21plusunique1} (D).  If we dig the bypass along
$A_+$ and add it back along $A_-$ as indicated, the result is
exactly  $M_2=M_1 \backslash (\delta \times I)$ with
$\Gamma_A=T2_1^+$ as pictured in figure \ref{T11T21plusunique1} (H).
Third, consider $M_3=M_2 \backslash (\delta_1 \times I)$ as pictured
in figure \ref{T11T21plusunique1}(E).  If we dig the bypass along
$A_+$ and add it back along $A_-$ as indicated, the result is
exactly $M_3=M_2 \backslash (\delta_1 \times I)$ with
$\Gamma_A=T2_1^+$ as pictured in figure \ref{T11T21plusunique1} (I).
Finally, consider $M_4=M_3 \backslash (\delta_2 \times I)$ as
pictured in figure \ref{T11T21plusST} (D).  If we dig the bypass and
add it back as indicated, the result is exactly  $M_4=M_3 \backslash
(\delta_2 \times I)$ with $\Gamma_A=T2_1^+$ with
$\#(\Gamma_{\partial M_4}) = 1$. Thus, we have that the unique
non-product $T1_1$ is equivalent to the unique non-product $T2_1^+$.

\begin{figure} [ht!]
\begin{center}
\includegraphics[bb=0 0 427 492, width=0.9\hsize]{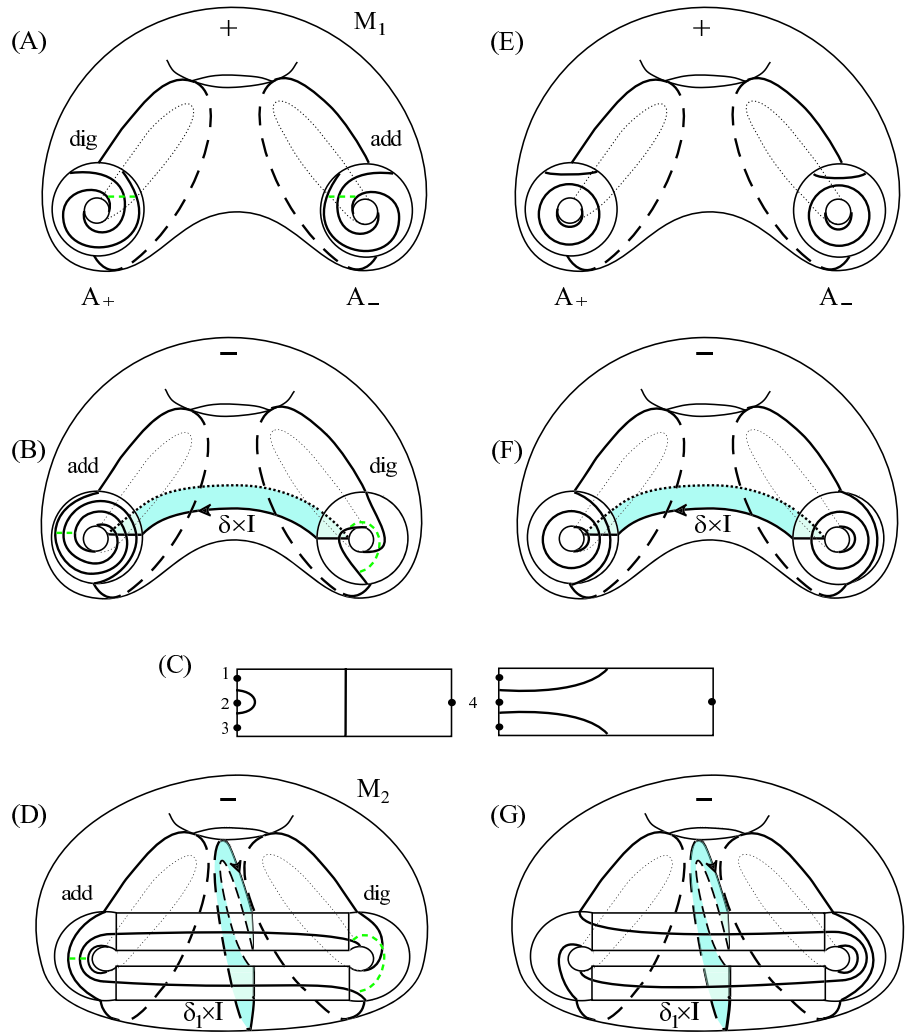}
\end{center}
\caption{\label{T11T21minusST1} Equivalence of $T1_1$ and $T2_1^-$
(1)}
\end{figure}

To obtain the equivalence $T2_1^- \cong T1_{-1}$, we use the same
convex decomposition as in the previous cases.  For the case $M_1$
with  $\Gamma_A =$  $T2_1^-$ and $T1_{-1}$, the convex decomposition
is given partly in figure \ref{T21minusT1-1unique1}. The rest of the
decomposition is as in figure \ref{T11T21minusST1}.

\begin{figure} [ht!]
\begin{center}
\includegraphics[bb=0 0 401 581,width=0.87\hsize]{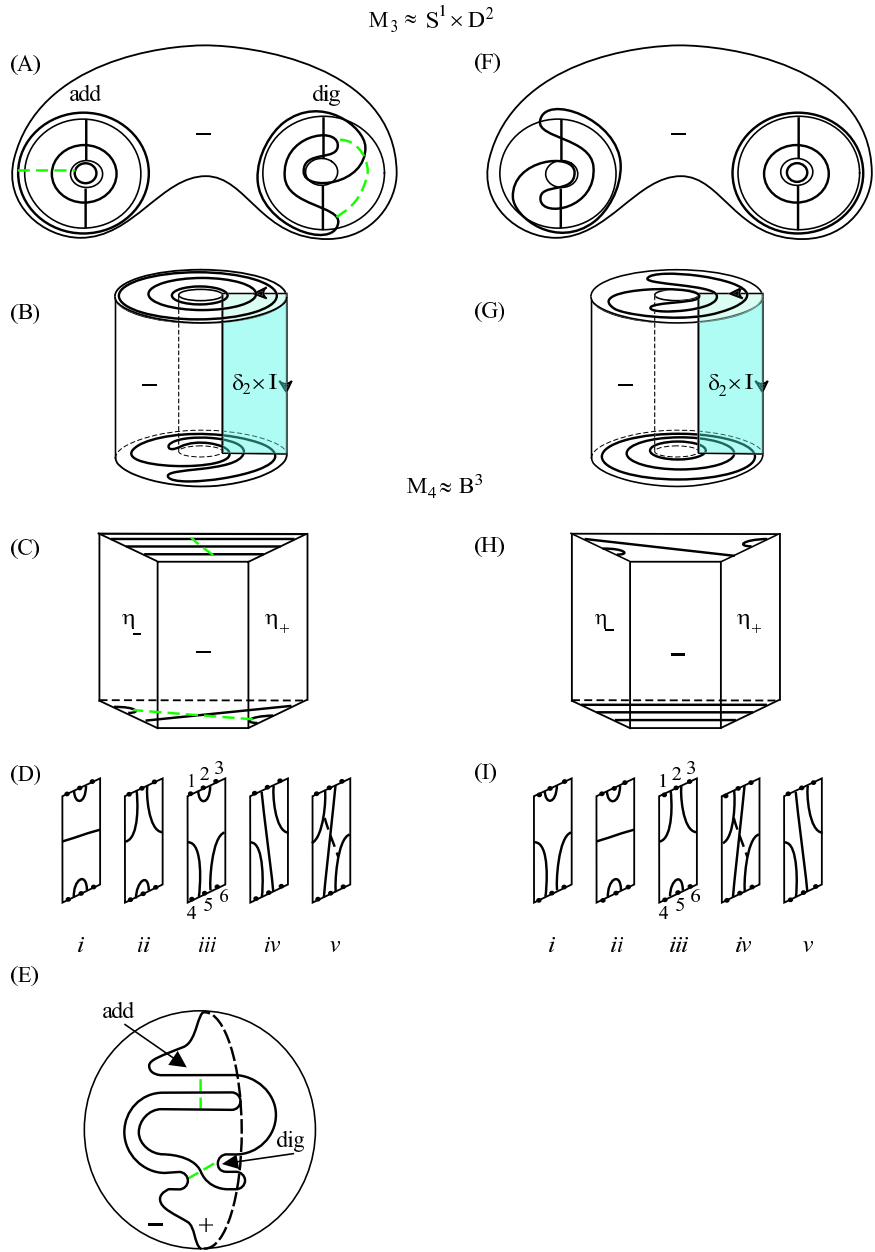}
\end{center}
\caption{\label{T11T21minusST2} Equivalence of $T1_1$ and $T2_1^-$
(2)}
\end{figure}

\begin{figure} [ht!]
\begin{center}
\includegraphics[bb=0 0 430 521,width=0.85\hsize]{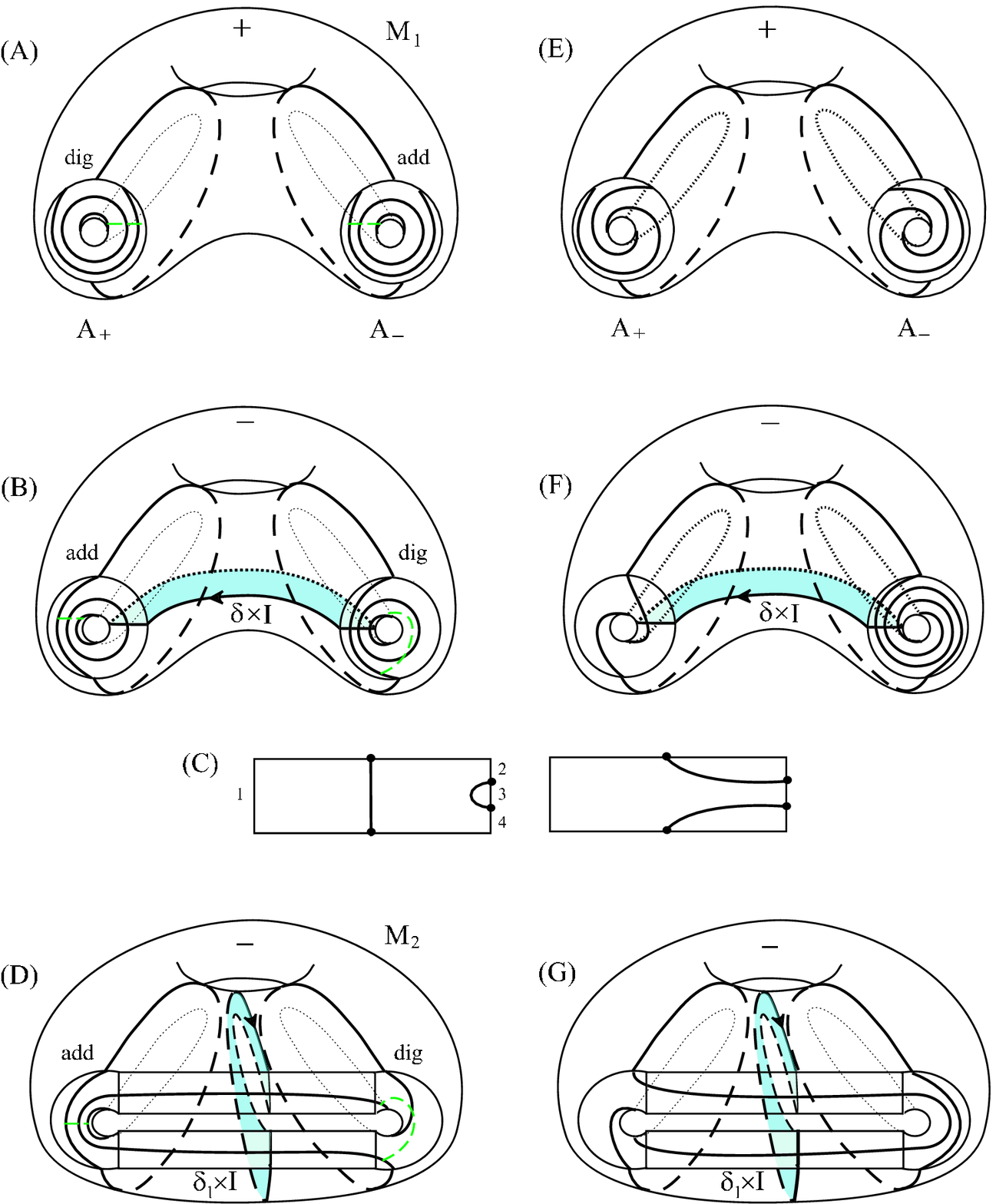}
\end{center}
\caption{\label{T21plusT1-1ST} Equivalence of $T2_1^+$ and
$T1_{-1}$}
\end{figure}

For the equivalences $T1_1 \cong T2_1^-$ and $T2_1^+ \cong
T1_{-1}$, we use a different convex decomposition to avoid cutting
through the proposed bypasses. The convex decomposition for the
first of these two transitions is given in figures
\ref{T11T21minusST1} and \ref{T11T21minusST2}. For the equivalence
$T2_1^+ \cong T1_{-1}$, the convex decomposition is given partly
in figure \ref{T21plusT1-1ST}.  The rest of the decomposition is
as in figure \ref{T11T21minusST2}.  \endproof

By our previous classification, we may conclude that the non-product
structure we have found is the unique non-product tight contact
structure on $M$.  Moreover, this decomposition shows us that the
non-product structure has a cut of type $T2_1^+$ and one of type
$T2_1^-$.  So, if either boundary component, $\delta \times \{0\}$
or $\delta \times \{1\}$, of our initial annular cutting surface is
the attaching arc of of a bypass, then our manifold contains all
possible complementary bypasses.  Thus, this attachment would result
in an overtwisted structure.

We can find a diffeomorphism of the surface, fixing the dividing
set, taking this $\delta$ to any other curve that can be the
attaching arc of a non-trivial bypass, and this diffeomorphism may
be extended to the product. Since there is unique non-product
structure, this must take the non-product structure $\xi$ to the
non-product structure, and, hence, all bypasses must exist along the
image of $\delta$. It follows that this manifold contains every
allowable bypass abutting the boundary.  This means that if
$(M,\xi)$ is contained in some tight $(M^\prime,\xi^\prime)$, then
for any convex surface with boundary $S \subset M^\prime$ such that
$\partial S \subset \partial M$ and $\# (\partial S \cap
\Gamma_{\partial M}) = 2$, the dividing set on $S$ cannot contain
any boundary-parallel dividing arcs. We may conclude that the
neighborhood of this cutting surface $\Sigma_2$ inside a tight
$(M^\prime, \xi^\prime)$ must be a product.  

\paragraph{\textbf{Acknowledgements}} I would like to thank
my doctoral advisor, Gordana Mati\'{c}, for her guidance and helpful
comments. I would also like to thank Joan Birman for inspiring me
and for providing indispensable advice and direction. My sincere
gratitude goes to Nancy Wrinkle for her unflagging support,
encouragement, and assistance throughout the entire process. I offer
my gratitude also to the referee, whose detailed and thoughtful
comments dramatically improved this paper. Finally, I would like to
thank Will Kazez, Ko Honda, Clint McCrory, John Hollingsworth as
well as numerous others for their significant contributions to my
education.

\Addresses\recd


\begin{thebibliography}

\bibitem{Benn}
\textbf{D Bennequin}, \emph{Entrelacements et \'equations de Pfaff},
  Ast\'{e}risque 107-108 (1983) 87--161 \MR{0753131}

\bibitem{colgl}
\textbf{V Colin}, \emph{Recollement de vari\'et\'es de contact tendues}, Bull.
  Soc. Math. France 127 (1999) 43--69
\MR{1700468} 

\bibitem{atoroidalC}
\textbf{V Colin}, \emph{Une infinit\'e de structures de contact tendues sur les
  vari\'et\'es toro{\"i}odales}, Comment. Math. Helv. 76 (2001) 353--372
\MR{1839351}

\bibitem{atoroidal2}
\textbf{V Colin}, \textbf{E Giroux}, \textbf{K Honda}, \emph{On the coarse
  classification of tight contact structures}, from: ``Topology and geometry of manifolds (Athens, GA, 2001)'', Proc. Sympos. Pure Math. 71, 
Amer. Math. Soc. (2003)
\MR{2024632}


\bibitem{EliashOTCS}
\textbf{Y Eliashberg}, \emph{Classification of overtwisted contact structures
  on 3-manifolds}, Invent. Math. 98 (1989) 623--637
\MR{1022310}

\bibitem{Eliash}
\textbf{Y Eliashberg}, \emph{Contact 3-manifolds twenty years since J.
  Martinet's work}, Ann. Inst. Fourier (Grenoble) 42 (1992) 165--192
\MR{1162559}

\bibitem{EtnyreLensSp}
\textbf{J Etnyre}, \emph{Tight contact structures on lens spaces}, Commun.
  Contemp. Math. 2 (2000) 559--577
\MR{1806947}

\bibitem{Nonexist}
\textbf{J Etnyre}, \textbf{K Honda}, \emph{On the nonexistence of tight contact
  structures}, Ann. of Math. 153 (2001) 749--766
\MR{1836287}

\bibitem{convexity}
\textbf{E Giroux}, \emph{Convexit\'e en topologie de contact}, Comment. Math.
  Helv. 66 (1991) 637--677
\MR{1129802}

\bibitem{GirLensSp}
\textbf{E Giroux}, \emph{Structures de contact en dimension trois et
  bifurcations des feuilletages de surfaces}, Invent. Math. 141 (2000) 615--689
\MR{1779622}

\bibitem{GirCircBund}
\textbf{E Giroux}, \emph{Structures de contact sur les vari\'et\'es fibr\'ees
  en cercles au-dessus d'une surface}, Comment. Math. Helv. 76 (2001) 218--262
\MR{1839364}

\bibitem{TCSI}
\textbf{K Honda}, \emph{On the classifiction of tight contact structures I},
\gtref4{2000}{11}{309}{368}
\MR{1786111}

\bibitem{TCSII}
\textbf{K Honda}, \emph{On the classifiction of tight contact structures II},
  J. Differential Geom. 55 (2000) 83--143
\MR{1849027}

\bibitem{gluing}
\textbf{K Honda}, \emph{Gluing tight contact structures}, Duke Math. J. 115
  (2002) 435--478
\MR{1940490}

\bibitem{TCSTF}
\textbf{K Honda}, \textbf{W Kazez}, \textbf{G Mati\'{c}}, \emph{Tight contact
  structures and taut foliations}, \gtref4{2000}7{219}{242}
\MR{1780749}

\bibitem{convexdecomp}
\textbf{K Honda}, \textbf{W Kazez}, \textbf{G Mati\'{c}}, \emph{Convex
  decomposition theory}, Internat. Math. Res. Notices  (2002) 55--88
\MR{1874319}

\bibitem{FHM}
\textbf{K Honda}, \textbf{W Kazez}, \textbf{G Mati\'{c}}, \emph{Tight contact
  structures on fibred hyperbolic 3-manifolds}, J. Differential Geom. 64 (2003)
  305--358
\MR{2029907}

\bibitem{haken}
\textbf{W Jaco}, \emph{Lectures on three-manifold topology}, CBMS Regional
  Conference Series in Mathematics 43, Amer. Math. Soc.
  Providence, RI (1980)
\MR{1874319}

\bibitem{KandaTor}
\textbf{Y Kanda}, \emph{The classification of tight contact structures on the
  3-torus}, Comm. in Anal. and Geom. 5 (1997) 413--438
\MR{1487723}

\bibitem{KandaLR}
\textbf{Y Kanda}, \emph{On the Thurson-Bennequin invariant of Legendrian knots
  and non exactness of Bennequin's inequality}, Invent. Math. 133 (1998)
  227--242
\MR{1632790}

\bibitem{ML}
\textbf{S Makar-Limanov}, \emph{Tight contact structures on solid tori}, Trans.
  Amer. Math. Soc. 350 (1998) 1045--1078
\MR{1401526}

\bibitem{Mart} \textbf{J Martinet}, \emph{Formes de contact sur les
vari\'et\'es de dimension 3}, from: ``Proceedings of Liverpool
Singularities Symposium, II (1969/1970)'', Lect. Notes in Math. 209,
Springer (1971) 142--163
\MR{0350771}


\end{thebibliography}
\end{document}